\documentclass{article}
\usepackage[numbers]{natbib}

\usepackage[preprint]{neurips_2024}

\usepackage[utf8]{inputenc} 
\usepackage[T1]{fontenc}    
\usepackage{hyperref}       
\usepackage{url}            
\usepackage{booktabs}       
\usepackage{amsfonts}       
\usepackage{nicefrac}       
\usepackage{microtype}      
\usepackage{xcolor}         

\usepackage{amssymb}   
\usepackage{amsthm}
\usepackage{graphicx}
\usepackage{stmaryrd}
\usepackage{nccmath}
\usepackage{geometry}
\usepackage{amsmath}
\newcommand{\R}{\mathbb{R}}
\newcommand{\C}{\mathbb{C}}
\newcommand{\N}{\mathbb{N}}
\newcommand{\E}{\mathbb{E}}

\usepackage{physics}

\DeclareMathOperator*{\argmin}{argmin}

\DeclareMathOperator*{\Diag}{Diag}

\newtheorem{thm}{Theorem}
\newtheorem{ppt}{Proposition}

\newtheorem{dft}{Definition}

\newtheorem{rmk}{Remark}

\title{WeSpeR: Computing non-linear shrinkage formulas for the weighted sample covariance}
\author{
  Benoît Oriol \\
  CEREMADE, Université Paris-Dauphine, PSL, Paris, France \\
Société Générale Corporate and Investment Banking, Puteaux, France \\
  \texttt{benoit.oriol@dauphine.eu} \\
}
\begin{document}

\maketitle

\begin{abstract}
We address the issue of computing the non-linear shrinkage formulas for the weighted sample covariance in high dimension. We use theoretical properties of the asymptotic sample spectrum in order to derive the \textit{WeSpeR} algorithm and significantly speed up non-linear shrinkage in dimension higher than $1000$. Empirical tests confirm the good properties of the \textit{WeSpeR} algorithm. We provide the implementation in PyTorch for it.
\end{abstract}

\section{Introduction and related work}
In Random Matrix Theory (RMT), the sample covariance spectrum has a non-random limit, denoted $F$, when the dimension grows linearly with the number of samples. In this regime, the sample spectral distribution does not converge to the population spectral distribution, but to a limit described by the Marcenko-Pastur equation \cite{Marcenko1967}. The joint work of Silverstein and Choi \cite{Silverstein1995c} gives important results on the asymptotic distribution $F$: $F$ has a density on $\R^*$, and its support can be computed with a simple procedure without any sampling of high dimensional sample covariance matrices. 

Those results were directly used in the design of the algorithm QuEST by Ledoit and Wolf \cite{Ledoit2015}, aiming at retrieving the population covariance spectrum and the sample covariance asymptotic density $F'$, and derive the non-linear shrinkage formulas theoretically given in \cite{Ledoit2009}.

Using a closed-form formula linking the Cauchy-Stieltjes transform of the empirical spectral density and the non-linear shrinkage formulas, analytical kernel estimators were derived to make the computation faster for high-dimensional data, for a reduced cost of performance \cite{Ledoit2019, Ledoit2020b}.

Those works focus on the sample covariance matrix, however in practice we often face \textbf{weighted} sample covariance matrices, in particular in multivariate time series analysis. Indeed, weighting schemes, such as the exponential weighted moving average (EWMA), are a model-free approach, and represent a widely used method to estimate statistics. They were used in covariance estimation for portfolio management in \cite{Pafka2004}, further studied for covariance filtering in \cite{Daly2010}, for financial spectrum estimation in \cite{Svensson2007}, and recently Tan et al.\cite{Tan2023} developed a NERCOME-like approach for EWMA sample covariance in a dynamic brain connectivity setting. Weighted sample covariance matrices are also directly observed in finance due to stochastic variance processes \cite{Bun2016}, and in wireless network they emerge in the model of Multiple Inputs Multiple Outputs (MIMO) systems \cite{Couillet2011}.

The theoretical extensions of the asymptotic results for the spectral distribution and the estimation of the weighted sample covariance matrices have been explored in various works. \cite{Zhang2007} establishes a Fundamental Equation analogous to the Marcenko-Pastur equation, proving that the spectrum converges to a deterministic distribution, denoted $F$, in high dimension. The analytical properties of $F$ and a theoretical framework for determining its support $S_F$ are developed in \cite{Couillet2015}. Additionally, asymptotic non-linear shrinkage formulas for estimating the covariance and precision matrices are presented in \cite{Oriol2025b}.

Those formulas are non-trivial to compute in practice as they depend on the unobserved population spectral distribution $H$ and are defined partially implicitly: the term $\check X$ is the limit on the real line of the solution of a fixed point equation in $\C_+ := \{z \in \C | \Im[z] > 0\}$. To address this issue, \cite{Oriol2025b} (Proposition 4) proposes a numerical procedure to estimate the population spectral distribution $H$ through a first-order minimization problem and a numerical procedure to compute the term $\check X$ (Theorem 2.2 and 2.3 \cite{Oriol2025b}). The minimization consists of sampling large weighted sample covariance matrices, computing their eigenvalues, and minimizing the distance between the observed and the sampled eigenvalues. Each step of the minimization problem has a complexity of $O(n^3)$ with an observed weighted sample covariance of size $(n,n)$.

While the implementation of this procedure is fairly simple, the cubic complexity in the dimension makes it unusable in practice for dimension $n$ exceeding a few thousands. We identify three regimes of dimension and propose a way of computing the non-linear shrinkage formulas in each case, leveraging analytical knowledge of the support of $F$ and stochastic Lanczos quadrature for the higher dimensional regimes:
\begin{itemize}
	\item[(LD)] the low dimensional case, when the cost of computing the eigenvalues of a $(n,n)$ Hermitian matrix at each minimization step is acceptable, roughly $n < 1000$,
	\item[(MD)] the medium dimensional case, when the cost of diagonalizing the observed weighted sample covariance of size $(n,n)$ once is acceptable, roughly $n < 20000$,
	\item[(HD)] the high dimensional case, when no routine of complexity $O(n^3)$ is acceptable, roughly $n > 20000$.
\end{itemize}

This work aims at providing an algorithm to compute asymptotic non-linear shrinkage in high dimension, with the implementation provided in this \href{https://www.github.com/nlcvbo/WeSpeR}{Github repository}, and a guide for the practitioner on which method to choose depending on its case. 

\section[Setting, formulas and (LD) procedure]{The setting, the formulas and the low dimensional procedure}
In this work, we consider studying the weighted sample covariance $S_n$, which is a random Hermitian matrix of size $(n,n)$ which can be decomposed as $S = \frac{1}{N} \sqrt{\Sigma_n} Z_n W_n Z_n^* \sqrt{\Sigma_n}$ where:
\begin{itemize}
	\item $\Sigma_n$ is a Hermitian nonnegative $(n,n)$ matrix: this is the population covariance matrix,
	\item $Z_n$ is a $(n,N)$ matrix with i.i.d. entries of mean zero and unit variance and bounded $12^{th}$ moment: this is the noise matrix,
	\item $W_n$ is a diagonal nonnegative $(N,N)$ matrix: this is the weight matrix.
\end{itemize}
We assume this structure in all that follows.

The asymptotic analysis of the spectrum of $S$ is conducted in the high dimensional regime: we assume that $n \rightarrow +\infty$ and $N := N(n)$ such that $n/N \rightarrow c > 0$, where usually $n$ is the dimension, $N$ the number of samples and $c$ is called the concentration ratio. 

A suitable object for the asymptotic analysis of the spectrum of a Hermitian matrix $M$ of size $(n,n)$ is its empirical spectral distribution (e.s.d.) denoted $F^M$ and defined as following: denoting $\lambda_1,...,\lambda_n$ the eigenvalues of $M$, we define $F^M = \frac{1}{n} \sum_{i=1}^n \mathbf{1}_{[\lambda_i,+\infty[}$, which is the cumulative distribution function of the spectrum of $M$. For readability, we will be using the notation used in \cite{Oriol2025b}: $H_n := F^{\Sigma_n}$, $D_n := F^{W_n}$ and $F_n := F^{S_n}$.

\subsection{Spectrum convergence and asymptotic non-linear shrinkage}

We give two central theoretical results that will be of importance in the design and the understanding of the numerical implementations: the Fundamental Equation from Theorem 1.2.1 \cite{Zhang2007} which establishes the convergence of $F^{S_n}$ to a non-random distribution $F$ and the non-linear shrinkage formulas for covariance estimation from Theorem 2.3 \cite{Oriol2025b} (we do not detail the formulas for the precision matrix for brevity, but the idea is similar). The two results are given with the notation of this work, and with a sufficient set of assumptions for our setting, for clarity, however not always optimal regarding the literature of each specific theorem.

\begin{thm}[Fundamental Equation, (from Theorem 1.2.1 \cite{Zhang2007})]\label{FE}
	Suppose:
	\begin{itemize}
		\item $H_n \underset{n \rightarrow +\infty}{\implies} H$ a.s. where $\implies$ denotes weak convergence, $H$ denotes a cumulative distribution function (c.d.f.),
		\item $D_n \underset{n \rightarrow +\infty}{\implies} D$ a.s., where $D$ denotes a c.d.f.,
		\item $\Sigma_n$, $Z_n$ and $W_n$ are mutually independent.
	\end{itemize}
	Then, $F_n \underset{n \rightarrow +\infty}{\implies} F$ a.s., where $F$ denotes a c.d.f., which is characterized by its Cauchy-Stieltjes transform $m(\cdot) := m_F(\cdot)$:
	\begin{equation}
	\begin{aligned}
		\forall z \in \C_+, m(z) = -\frac{1}{z} \int \frac{1}{\tau X(z) + 1}dH(\tau),
	\end{aligned}
	\end{equation} 
	where $X(z)$ is the unique solution in $\C_+$ of the fixed point equation:
	\begin{equation}
	\begin{aligned}
		X(z) = -\int \frac{\delta}{z - \delta c \int \frac{\tau}{\tau X(z) + 1}dH(\tau)}dD(\delta).
	\end{aligned}
	\end{equation} 
\end{thm}

We introduced the central asymptotic objects of non-linear shrinkage, $H$, $D$, $c$, $F$, and $X(\cdot)$, and we can state the formulas for asymptotic non-linear shrinkage. We will remain intentionally vague in the formulation of "asymptotic non-linear shrinkage formulas for covariance estimation" for brevity, the concepts are developed and explained initially in \cite{Ledoit2009} for the sample covariance and then in \cite{Oriol2025b} for the weighted sample covariance.

\begin{thm}[Non-linear shrinkage formulas for covariance estimation (from Theorem 2.3, \cite{Oriol2025b})]
	Suppose:
	\begin{itemize}
		\item $H_n \underset{n \rightarrow +\infty}{\implies} H$ a.s., where $H$ denotes a c.d.f., and the support of $H_n$, denoted $S_{H_n}$ is included in $[h_1,h_2]$ where $0<h_1\leq h_2 < \infty$,
		\item $D_n \underset{n \rightarrow +\infty}{\implies} D$ a.s., where $D$ denotes a c.d.f., and $S_{D_n} \subset [d_1, d_2]$ where $0 < d_1 \leq d_2 < \infty$,
		\item $\Sigma_n$, $Z_n$ and $W_n$ are mutually independent.
	\end{itemize}
	Then, for $\lambda \in \R$, the asymptotic non-linear shrinkage formula at $\lambda$ for covariance estimation, denoted $h(\lambda)$, is:
	\begin{equation}\label{nl_formulas}
	\begin{aligned}
		& h(\lambda) = \frac{\int \frac{\tau^2}{|\tau \check X(\lambda) + 1|^2}dH(\tau)}{\int \frac{\tau}{|\tau \check X(\lambda) + 1|^2}dH(\tau)}, \text{ if } \lambda \neq 0, \\
		& h(\lambda) = \frac{1}{(c-1) \check X(0)}, \text{ if } \lambda = 0 \text{ and } c > 1, \\
		& h(\lambda) = 0, \text{ otherwise,}
	\end{aligned}
	\end{equation} 
	where $\check X(\lambda) = \lim_{\eta \rightarrow 0^+} X(\lambda + i\eta)$.
\end{thm}

This formula means that for an asymptotic optimal estimation of the covariance in the class of rotation-invariant estimators (the estimators that have the same eigenvectors as the weighted sample covariance $S_n$), we need to replace the sample eigenvalues $\lambda_i$ by $h(\lambda_i)$ in $S_n$ (and keep the same eigenvectors). In matrix form, if we observe $S_n = U_n \Lambda_n U_n^*$ under its diagonalized form, the asymptotic optimal covariance estimator is $S_n^* := U_n h(\Lambda_n) U_n^*$, where $h(\Lambda_n) = \Diag\left(h(\lambda_1),...,h(\lambda_n)\right)$.

Similarly, for the estimation of the precision matrix $\Gamma_n := \Sigma_n^{-1}$, the asymptotic optimal precision estimator is $P_n^* := U_n t(\Lambda_n) U_n^*$, where the function $t$ is defined in Theorem 2.5 \cite{Oriol2025b}, and can be computed with $H$, $\check X$ and $c$, as $h(\cdot)$.

\subsection{Low-dimensional (LD) procedure for computing the shrinkage formulas}
As a benchmark and introduction to the numerical implementation of non-linear shrinkage for weighted sample covariance, we remind the procedure given in Section 5 \cite{Oriol2025b}.

The idea is to estimate the population spectrum distribution $H$ as a mixture of Dirac from the observed weighted sample eigenvalues $(\lambda_i)_{i=1}^n$, and using this distribution $H$ to compute $\check X$, and the shrinkage formulas.

\begin{itemize}
	\item[1-] As input, we take the observed sample spectrum distribution $F_n = \frac{1}{n} \sum_{i=1}^n 1_{[\lambda_i,\infty[}$ and the weight matrix $W$.
	\item[2-] Find the estimated population spectrum $\hat H(\tau) = \frac{1}{n}\sum_{i=1}^n 1_{[\tau_i,+\infty[}$ where $\tau$ solves:
		\begin{equation}\label{autoquest}
		\begin{aligned}
			\min_{\tau \in \R^n}\E_{Z}\left[\left\lVert \tilde F_n(Z, \hat H) - F_n \right\rVert_{\mathcal{W},2}^2 \right]
		\end{aligned}
		\end{equation} 
	where $\lVert \cdot \rVert_{\mathcal{W},2}$ is the $2$-Wasserstein norm and $\tilde F_n(Z) := F^{\frac{1}{N}\sqrt{T}ZWZ^*\sqrt{T}}$ with:
	\begin{itemize}
		\item $T = \text{Diag}\left((\tau_i)_{i=1}^n\right)$, 
		\item $Z$ of size $(n,N)$ with i.i.d. $Z_{ij} \sim \mathcal{N}(0,1)$.
	\end{itemize}
	We use automatic differentiation along with Adam optimizer to solve it.
	\item[3-] Compute $\check X(\lambda_i)$ with Proposition 4 \cite{Oriol2025b}, solving if $\lambda_i \in S_F$: 
		\begin{equation}
		\begin{aligned}
			\check X(\lambda_i) = \argmin_{X \in \C_+} \left| X + \int \frac{\delta}{\lambda_i - \delta c \int \frac{\tau}{\tau X + 1}d\hat H(\tau)}dD(\delta) \right|^2
		\end{aligned}
		\end{equation} 
		Remark that for $n$ large enough, $\lambda_i \in S_F$ with probability one \cite{Paul2009}, but if $\lambda_i \notin S_F$, the minimization has no solution in $\C_+$, in this case define $\check X(\lambda_i) = \Re[X(\lambda_i + i\epsilon)]$ for a numerically small $\epsilon > 0$.
	\item[4-] Compute $h(\lambda_i)$ using Equation \ref{nl_formulas} with $\hat H$ and $\check X(\lambda_i)$.
\end{itemize}

\begin{rmk}[Noise sampling]
We note that the expectation in \eqref{autoquest} can be computed under any centered and standardized distribution, irrespectively of the noise of the observed phenomenon due to the universality the asymptotic spectrum.
\end{rmk}

\begin{figure}[]
\centering
\includegraphics[width=0.7\linewidth]{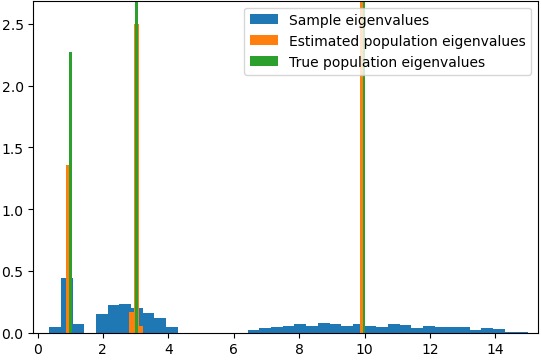}
\includegraphics[width=0.7\linewidth]{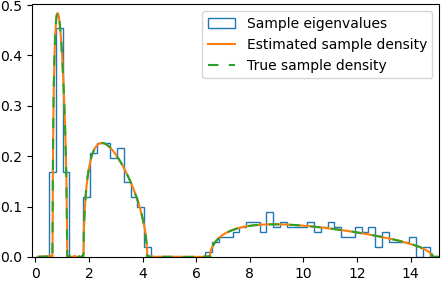}
\caption[Histograms and densities resulted by (LD) procedure.]{(Top) Histograms of sample eigenvalues, estimated population eigenvalues $\hat H$ with the LD procedure, and true population eigenvalues $H$. (Bottom) Estimated and theoretical sample density and sample eigenvalues' histogram.}
\label{fig:autoquest2}
\end{figure}

An experimental result is shown in Figure \ref{fig:autoquest2} to illustrate how the algorithm works. Instead of plotting $h(\lambda_i)$, we show for clarity the estimated asymptotic density $F': \lambda \in \R^* \mapsto \frac{1}{\pi} \int \frac{\tau \Im[\check X(\lambda)]}{\lambda |\tau \check X(\lambda) + 1|^2}d\hat H(\tau)$ along with the histogram of observed eigenvalues $\lambda_i$. We used $H = \frac{1}{5} \mathbf{1}_{[1,\infty[} + \frac{2}{5} \mathbf{1}_{[3,\infty[} + \frac{2}{5} \mathbf{1}_{[10,\infty[}$, $D$ exponentially weighted with $\alpha = 1$, $c=0.1$, and $Z_{ij} \sim  \mathcal{N}(0,1)$.

While this algorithm is very efficient in practice for low dimension, it struggles in higher dimension for two reasons. In fact, each gradient descent step in step (2) involves:
\begin{itemize}
	\item sampling a random noise matrix $Z$ of size $(n,N)$, which causes a memory and time complexity of order $O(n^2/c)$ which is problematic even for quite low dimension $n \approx 1000$ if $c$ is low. Remark that if the weight distribution $D$ has a large variance, typically EWMA distribution with large decay parameter $\alpha$, non-linear shrinkage still has a significant impact in practice even with very low $c$,
	\item computing the eigenvalues of $\frac{1}{N}\sqrt{T}ZWZ^*\sqrt{T}$ while building the backward graph is of time complexity $O(n^3)$, and it scales badly with dimension $n > 1000$.
\end{itemize}

\section[The medium dimensional (MD) algorithm]{The medium dimensional (MD) algorithm: identification of $S_F$ and inversion of the Fundamental Equation with \textit{WeSpeR}}
The goal of this section in to provide an algorithm that gets rid of the $O(n^3)$ complexity at each step of the minimization to retrieve $H$ as a mixture of Dirac. The general idea of the algorithm remains the same as the LD procedure, we will only, but deeply, change step (2). While, in (LD), $\tilde F_n := \tilde F_n(Z, \hat H)$ was obtained through sampling, this time we are computing $\tilde F := \tilde F(c, \hat H)$ solving the Fundamental Equation described in Theorem \ref{FE}.

The idea to "invert" the Fundamental Equation in order to retrieve the population spectrum $H$ was introduced in the QuEST algorithm \cite{Ledoit2015} for the unweighted sample covariance. We generalize the approach to the weighted sample covariance in this algorithm, denoted \textit{WeSpeR} for \textit{W}eighted sample covariance \textit{Spe}ctrum \textit{R}etrieval algorithm.

As in QuEST, to enable a first-order minimization scheme, we compute analytically the derivatives. Because we want to rely on non-commercial software and optimization routines, our implementation architecture differs from QuEST, originally in Matlab with SNOPT optimization routine. We implemented the algorithm in Python, as a PyTorch module with its own backward function analytically computing the derivatives, and optimizing through PyTorch optimizers (by default Adam optimizer). The algorithm can then be integrated in larger and more complex optimization problems, using automatic differentiation and the PyTorch engine.

The algorithm to compute $\tilde F(c, \hat H)$ functions as following, supposing a population spectrum $\hat H = \sum_{i=1}^n \mathbf{1}_{[\tau_i,+\infty[}$:
\begin{enumerate}
	\item Find the support $S_{\tilde F}$ of $\tilde F$. We note $S_{\tilde F} = \cup_{i=1}^\nu [s_{2i-1},s_{2i}]$.
	\item Choose a grid $\xi_i^j$ that covers the support $S_{\tilde F}$ following an arcsine distribution on each $[s_{2i-1},s_{2i}]$.
	\item Solve the Fundamental Equation, \textit{i.e.}, compute $\check X(\xi_i^j)$ in this grid.
	\item Compute the spectral density $\tilde F'(\xi_i^j) = \pi^{-1}\Im[\check m(\xi_i^j)]$, and integrate it to obtain the empirical $\tilde F$.
	\item Compute the loss $\ell(\tilde F(c, \hat H), F)$ and its derivatives. Update $\hat H$ with an order one descent algorithm and go back to step (1).
\end{enumerate}

Several issues have to be addressed to make this procedure possible to implement in practice:
\begin{itemize}
	\item in step (1), we need to find the support $S_{\tilde F}$ of the asymptotic distribution $\tilde F$, knowing $\hat H$, $c$, and the weight distribution $D$. This is essential to define a discretization grid. We do not want to sample and diagonalize a large matrix $(n,n)$ to estimate the support $S_{\tilde F}$, so we are using an analytical procedure based on a theoretical result in \cite{Couillet2015}, developed in the next section, to determine $S_{\tilde F}$,
	\item use a suitable loss $\ell(\tilde F(c, \hat H), F)$,
	\item compute analytically the necessary derivatives in the procedure to enable the minimization over $\tau$.
\end{itemize}

We present the general theoretical method to find the asymptotic support $S_{F}$ in function of $H$, $c$, and $D$. Then, we precisely detail the algorithm and its derivatives step by step.

\subsection{Identification of $S_F$}\label{secSF}
The purpose of this section is to find the support of $F$, denoted $S_F$, using a theoretical result developed and proved in \cite{Couillet2015}. The idea is to use one or several well-chosen real functions, easy to compute, and deduce the border of $S_F$ from the zero's of their derivatives. This method do not rely on sampling any weighted sample matrix and can detect even very small spectral gaps.

The theorems are given using the setting and notation of this work for coherence, and in a way to make the implementation straightforward.

Let us start where the support of the weight distribution, denoted $S_D$, is convex. 

\subsubsection{Identification of the support of $F$ in function of $H$ when $S_D$ is convex}
This case is very similar to the unweighted scenario, studied in \cite{Silverstein1995c}. Indeed, the number of spectral gaps in $S_F$ is bounded by the number of gaps in $S_H$, and one function is enough to detect all of them. The function we are interested in to determine $S_F$ is $x_F$, defined below. 
\begin{dft}
	Suppose $S_D$ is convex, \textit{i.e.} $S_D$ is of the form $[d_1,d_2]$. We define: 
	\begin{equation*}
	\begin{aligned}
		m_{LD}: x \in S_D^c \mapsto \int \frac{\delta}{\delta - x}dD(\delta) \in \R^*.
	\end{aligned}
	\end{equation*}
	Notice that $m_{LD}$ is invertible.
	Moreover, we define, with $B = \{y \in \R, y \neq 0, -\frac{1}{y} \in S_H^c \}$: 
	\begin{equation*}\label{}
	\begin{aligned}
		&x_F: X \in B \mapsto 
			\begin{cases}
				\frac{h(X)}{X} m_{LD}^{-1} \left(h(X) \right) & \text{, if }  h(X) \neq 0,\\
				- \frac{1}{X}\int \delta dD(\delta) & \text{, otherwise.}
			\end{cases}\\
		&\text{with } h: X \in B \mapsto cX  \int \frac{\tau}{\tau X +1}dH(\tau).
	\end{aligned}
	\end{equation*}
\end{dft}

\begin{ppt}\label{pptc1}
	Suppose $S_D$ is convex. Then $x_F \in \mathcal{C}^1(B,\R)$.
\end{ppt}

The practical result in the following theorem links $S_F$ with the derivative of $x_F$. 
\begin{thm}[Spectrum identification, from Section 3.3 \cite{Couillet2015}]\label{spectrum1}
	Suppose:
	\begin{itemize}
		\item $H$ is a c.d.f. $\R_+$,
		\item $D$ is a c.d.f. on $\R_+$, and $S_D$ is convex, \textit{i.e.} $S_{D} = [d_1,d_2]$,
		\item $F$ solves the Fundamental Equation \ref{FE} with $H$, $D$ and $c \in \R^*$.
	\end{itemize} 
	Then, $x \in S_{F}^c \iff X \in B \text{ and } x_F'(X) > 0$, where $X$ and $x$ are linked respectively by $X = \check X(x)$ and $x = x_F(X)$.
\end{thm}

The procedure to use this theorem is simple:
\begin{itemize}
	\item find the open intervals $(]a_i, b_i[)_i$ of $B$, $a_i < b_i$,  where $\forall i, \forall X \in ]a_i, b_i[, x_F'(X) > 0$,
	\item then, $S_F^c = \underset{i}\cup ]x_F(a_i), x_F(b_i)[$.
\end{itemize}
Numerically, we only need to find the zeros of $x_F'$ and compute the value of $x_F$ at those points.

\begin{figure}[]
\centering
\includegraphics[width=0.7\linewidth]{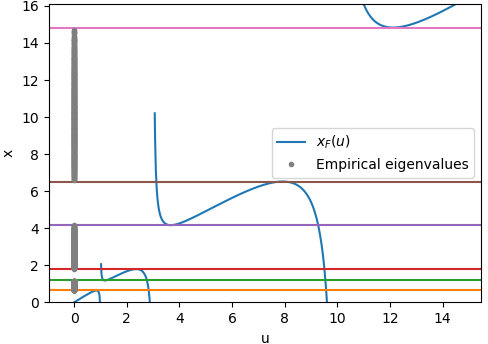}
\caption[Support identification procedure using $H$ mixture of $3$ Dirac, $D$ uniform.]{$x_F(-1/u)$ for $u \in S_H^c$, using $H = \frac{1}{5} \mathbf{1}_{[1,\infty[} + \frac{2}{5} \mathbf{1}_{[3,\infty[} + \frac{2}{5} \mathbf{1}_{[10,\infty[}$, $D$ uniform between $1/2$ and $3/2$, $c = 0.1$, and empirical eigenvalues sampled with $n=1000$, Gaussian noise.}
\label{fig:sep_1dirac}
\end{figure}

An example is given in Figure \ref{fig:sep_1dirac}, and illustrate the use of the theorem, and the precision of the prediction. Additional figures and experiments for diverse weight distributions are detailed in the Appendix.

\subsubsection{Identification of $S_F$ in function of $H$ when $S_D$ is a finite union of intervals}
In the more general case where $S_D$ is a finite union of intervals, the method can be extended. This case is useful in the case of $D$ being a finite mixture of Dirac for example. This method requires as many functions $x_F^{(k)}$ as there are disjoint intervals in $S_D$.

\begin{dft}
	 Suppose $S_D$ finite union of $M \in \N^*$ intervals, \textit{i.e.} there exists $\delta_1^{(1)} \leq\delta_2^{(1)} <...<\delta_1^{(M)} \leq \delta_2^{(M)}$ such that $S_D = \cup_{k=1}^M [\delta_1^{(k)}, \delta_2^{(k)}]$. We define for $k \in \llbracket 1, M-1 \rrbracket$: 
	\begin{equation*}\label{}
	\begin{aligned}
		& m_{LD}^{(k)}: x \in ]\delta_2^{(k)},\delta_1^{(k+1)}[ \mapsto \int \frac{\delta}{\delta - x}dD(\delta),
	\end{aligned}
	\end{equation*} 
	and,
	\begin{equation*}\label{}
	\begin{aligned}
		& m_{LD}^{(M)}: x \in \R \backslash [\delta_1^{(1)},\delta_2^{(M)}] \mapsto \int \frac{\delta}{\delta - x}dD(\delta).
	\end{aligned}
	\end{equation*} 
	For all $k \in \llbracket 1, M \rrbracket$, $m_{LD}^{(k)}$ is invertible.
	Moreover, we define for $k \in \llbracket 1, M-1 \rrbracket$, with $B = \{y \in \R, y \neq 0, -\frac{1}{y} \in S_H^c \}$: 
	\begin{equation*}\label{}
	\begin{aligned}
		x_F^{(k)}: X \in B \mapsto \frac{h(X)}{X} \left(m_{LD}^{(k)}\right)^{-1} \left(h(X) \right),
	\end{aligned}
	\end{equation*} 
	and we define $x_F^{(M)}: X \rightarrow \R$ as, for $X \in B$:
	\begin{equation*}\label{}
	\begin{aligned}
		x_F^{(M)}(X) = 
			\begin{cases}
				\frac{h(X)}{X} \left(m_{LD}^{(M)}\right)^{-1} \left(h(X) \right) & \text{, if }  h(X) \neq 0,\\
				- \frac{1}{X}\int \delta dD(\delta) & \text{, otherwise.}
			\end{cases}
	\end{aligned}
	\end{equation*} 
\end{dft}

As Proposition \ref{pptc1}, we have the following result.
\begin{ppt}\label{pptc12}
	 Suppose $S_D$ is a finite union of $M \in \N^*$ intervals, then $\forall k \in \llbracket 1, M \rrbracket, x_F^{(k)} \in \mathcal{C}^1(B,\R)$.
\end{ppt}

In this scenario, the theorem linking $S_F$ to the $(x_F^{(k)})_k$ is similar to Theorem \ref{spectrum1}: if any of the $x_F^{(k)}$ is increasing on an interval $I \subset \R$, then $x_F^{(k)}(I) \subset S_F^c$. This is formally stated in the following theorem. 
\begin{thm}[Spectrum identification, from Section 3.3 \cite{Couillet2015}]\label{spectrum2}
	Suppose:
	\begin{itemize}
		\item $H$ is a c.d.f. function on $\R_+$,
		\item $D$ is a c.d.f. function on $\R_+$, and $S_D = \sum_{i=1}^M [d_{2i-1},d_{2i}]$ with $d_1 \leq d_2 < ... < d_{2M-1} \leq d_{2M}$,
		\item $F$ solves the Fundamental Equation \ref{FE} with $H$, $D$ and $c \in \R^*$.
	\end{itemize}
	Let $x \in \R^*$. Then:
	\begin{equation*}\label{}
	\begin{aligned}
		x \in S_{ F}^c \iff \exists k \in \llbracket 1,M \rrbracket, \exists X \in B,  \left(x_F^{(k)}\right)'(X) > 0 \text{ and } \left(x_F^{(k)}\right)(X) = x.
	\end{aligned}
	\end{equation*}
\end{thm}

This result highlights a phenomenon of higher spectral separation than the previous section where $S_D$ was made of only one interval. Let us look at the same example of distribution $H$ we studied previously: $H$ being a mixture of 3 Dirac in $1$, $3$, and $10$ with respectively weights $0.2$, $0.4$ and $0.4$. With $S_D$ made of one interval, $S_F$ is a union of at most $3$ distinct intervals. Now, the situation is different, each separation in $S_D$ can lead to a spectral separation in the empirical spectrum. 

\begin{figure}[]
\centering
\includegraphics[width=0.39\linewidth]{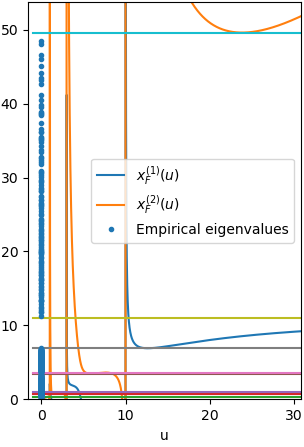}
\includegraphics[width=0.4\linewidth]{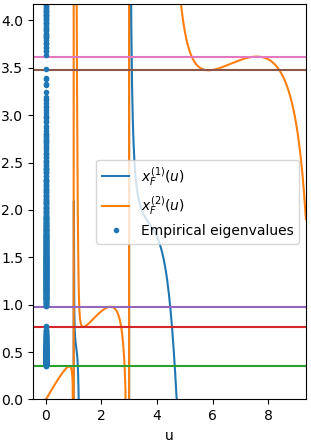}
\caption[Support identification procedure using $H$ mixture of $3$ Dirac, $D$ mixture of $2$ Dirac.]{$u \mapsto x_F^{(k)}(-1/u)$ for $u \in S_H^c$, using $H = \frac{1}{5} \mathbf{1}_{[1,\infty[} + \frac{2}{5} \mathbf{1}_{[3,\infty[} + \frac{2}{5} \mathbf{1}_{[10,\infty[}$, $D = \left(1-\frac{1}{80}\right)1_{[\frac{1}{2},\infty[} + \frac{1}{80}1_{[\frac{81}{2},\infty[}$, $c = 0.1$, and empirical eigenvalues sampled with $n=1000$. (Left) Whole spectrum. (Right) Zoom for small eigenvalues.}
\label{fig:sep_2Dirac}
\end{figure}

We show this specific behavior with weights following a mixture of two Dirac in Figure \ref{fig:sep_2Dirac} where $\mathbf{3}$ spectral gaps are visible in $S_F$ while $S_H$ has only $\mathbf{2}$ gaps. The case with $N$ Dirac is discussed in the Appendix, along with a way to compute it efficiently.

\begin{rmk}
	In practice, we use the increasing change of variable $u = -1/x$, and find the zeros of the derivatives of the functions $y_F^{(k)}: u \in S_H^c \backslash \{ 0 \} \mapsto x_F^{(k)}(-1/u)$.
\end{rmk}

\begin{rmk}\label{KD}
	It can be computationally demanding to compute each of the $x_F^{(k)'}$ and find its zeros if $M$ is large. A general heuristic we observed experimentally is that studying $x_F^{(M)'}$ is enough to find out most of the gaps if $S_D$ has no "large" gaps. We can see this heuristic as a "convex support" approximation of $S_D$.
\end{rmk}

The details of implementation of the routine to find the support $S_F$ is detailed in the Appendix for several weight distributions.

\subsection{\textit{WeSpeR}: retrieving $H$ without sampling}
Before detailing every step of the algorithm, we recall the sketch of it, supposing a population spectrum $\hat H = \sum_{i=1}^n \mathbf{1}_{[\tau_i,+\infty[}$:
\begin{enumerate}
	\item Find the support $S_{\tilde F}$ of $\tilde F$ using the numerical procedure developed previously. We note $S_{\tilde F} = \cup_{i=1}^\nu [s_{2i-1},s_{2i}]$.
	\item Choose a grid $\xi_i^j$ that covers the support $S_{\tilde F}$ following an arcsine distribution on each $[s_{2i-1},s_{2i}]$.
	\item Solve the Fundamental Equation, \textit{i.e.}, compute $\check X(\xi_i^j)$ in this grid.
	\item Compute the spectral density $\tilde F'(\xi_i^j) = \pi^{-1}\Im[\check m(\xi_i^j)]$, and integrate it to obtain the empirical $\tilde F$.
	\item Compute the loss $\ell(\tilde F, F)$ and its derivatives. Update $\hat H$ with an order one descent algorithm and go back to step (1).
\end{enumerate}

The main challenge of the implementation is to correctly compute the analytical derivatives with regard to the inputs $(\tau_i)_{i=1}^n$. Doing that, we can implement \textit{WeSpeR} as a PyTorch module, and integrate it in larger functions or any other differentiable loss effortlessly. We are computing the derivatives step-by-step. For generality purpose, we consider $\hat H = \sum_{k=1}^n w_k \mathbf{1}_{[\tau_k,+\infty[}$ with $w_k \geq 0$ and $\sum w_k = 1$.

\subsubsection{Step (1)}
Using the numerical procedure presented in the section \ref{secSF}, we can compute $s_1 < ... < s_{2\nu}$ such that $S_{\tilde F} = \cup_{i=1}^\nu [s_{2i-1},s_{2i}]$. Following the procedure to find them, we denote by $M$ the number of disjoint closed intervals forming $S_D$, then for each $i \in \llbracket 1, 2\nu \rrbracket$, there exists $m \in \llbracket 1, M \rrbracket$ and $u_i \in S_H^c$ such that $s_i = y_{\tilde F}^{(m)}(u_i)$ and ${y_{\tilde F}^{(m)}}'(u_i) = 0$.

Let $k \in \llbracket 1, n \rrbracket$. We want to compute the derivative ${\partial s_i}/{\partial \tau_k}$. We denote:
\begin{equation}\label{}
\begin{aligned}
	&y_{\tilde F}^{(m)}: (u; \tau_1,...,\tau_p) \mapsto x_{\tilde F}^{(m)}(-1/u; \tau_1,...,\tau_p) = -ut(u) m_{LD}^{(m)-1}(t(u)),\\
	&\text{with }t: u \in S_H^c \mapsto c\int \frac{\tau}{\tau-u}d\hat H(\tau).
\end{aligned}
\end{equation} 
We have:
\begin{equation}\label{}
\begin{aligned}
	&\frac{\partial s_i}{\partial \tau_k} = \underbrace{\frac{\partial y_{\tilde F}^{(m)}}{\partial u}(u_i)}_{= 0}\frac{\partial u_i}{\partial \tau_k} + \frac{\partial y_{\tilde F}^{(m)}}{\partial \tau_k} = \frac{\partial y_{\tilde F}^{(m)}}{\partial \tau_k}\\
	&\frac{\partial s_i}{\partial \tau_k} =  \frac{cu_i^2w_k}{(\tau_k - u_i)^2} \left(m_{LD}^{(m)-1}(t(u_i)) + \frac{t(u_i)}{m_{LD}'\left(m_{LD}^{(m)-1}(t(u_i))\right)} \right),
\end{aligned}
\end{equation} 
with $\forall x \in S_D^c, m_{LD}'(x) = \int\frac{\delta}{(\delta - x)^2}dD(\delta)$.

\subsubsection{Step (2)}
We have $S_{\tilde F} = \cup_{i=1}^\nu [s_{2i-1},s_{2i}]$. For $i \in \llbracket 1, \nu \rrbracket$, we extend the procedure used in QuEST \cite{Ledoit2015}: the interval $[s_{2i-1},s_{2i}]$ contains a proportion $\tilde \omega_i \in [0,1]$ of the sample eigenvalues \cite{Bai2004}, that can be estimated through $\hat H$ and $D$, and will be useful in the definition of the grid. We denote $S_{\hat H} = \cup_{\ell=1}^L [h_{2\ell-1},h_{2\ell}]$ and $S_D = \cup_{m=1}^M [d_{2m-1},h_{2m}]$. As in the previous step, there exists $\ell_1, \ell_2 \in \llbracket 0, L \rrbracket$, $m_1, m_2 \in \llbracket 0, M \rrbracket$ and $u_{2i-1} \in ]h_{2\ell_1},h_{2\ell_1+1}[, u_{2i} \in ]h_{2\ell_2},h_{2\ell_2+1}[$, with the convention $h_0 = -\infty, h_{2L+1} = +\infty, d_0 = -\infty, d_{2M+1} = +\infty$, such that $s_{2i-1} = y_{\tilde F}^{(m_1)}(u_{2i-1})$, and $s_{2i} = y_{\tilde F}^{(m_2)}(u_{2i})$. We define $\tilde \omega_i$ as following:
\begin{equation}\label{}
\begin{aligned}
	\tilde \omega_i = \sum_{\ell=\ell_1+1}^{\ell_2}\hat H([h_{2\ell-1},h_{2\ell}]) \sum_{j=m_1+1}^{m_2} D([d_{2k-1},d_{2k}]) .
\end{aligned}
\end{equation} 

As a hyperparameter, we define $\kappa \in \N^*$ the number of points in the discretization grid. We define then $\omega_i = \lceil \kappa \tilde \omega_i \rceil$: there will be $\omega_i+2$ points of discretization in the interval $[s_{2i-1},s_{2i}]$. We define the discretization grid $(\xi_i^j)_{j=0}^{\omega_i+1}$ as following:
\begin{equation}\label{}
\begin{aligned}
	\xi_i^j = s_{2i-1} + (s_{2i} - s_{2i-1}) \sin^2\left[\frac{\pi j}{2(\omega_i + 1)}\right].
\end{aligned}
\end{equation} 

There are two ideas behind this definition:
\begin{itemize}
	\item through $\omega_i$, we use more discretization points on the intervals where there are a lot of sample eigenvalues;
	\item through the arcsine distribution of the $\xi_i^j$ in $[s_{2i-1},s_{2i}]$, we put more density near the border of the interval, at an inverse square-root rate: it follows the idea that the asymptotic density $\tilde F$ has a square-root behavior near its borders \cite{Couillet2015}.
\end{itemize}

We have the following derivatives:
\begin{equation}\label{}
\begin{aligned}
	\frac{\partial \xi_i^j}{\partial \tau_k} = \left(1 - \sin^2\left[\frac{\pi j}{2(\omega_i + 1)}\right] \right) \frac{\partial s_{2i-1}}{\partial \tau_k} + \sin^2\left[\frac{\pi j}{2(\omega_i + 1)}\right]\frac{\partial s_{2i}}{\partial \tau_k}.
\end{aligned}
\end{equation} 

\subsubsection{Step (3)}
For $i \in \llbracket 1, \nu \rrbracket$, $j \in \llbracket 1, \omega_i \rrbracket$, we have that $\check X_i^j := \check X(\xi_i^j)$ is the unique solution in $\C_+$ of the equation $\Gamma_i^j(X; \tau_1,...,\tau_n) = 0$ where:
\begin{equation}\label{}
\begin{aligned}
	\Gamma_i^j(X; \tau_1,...,\tau_n) = X + \int \frac{\delta}{\xi_i^j(\tau_1,...,\tau_n) - \delta c \int \frac{\tau}{\tau X + 1}dH(\tau)}dD(\delta).
\end{aligned}
\end{equation} 
In order to compute the derivatives, we define:
\begin{equation}\label{}
\begin{aligned}
	&\hat \Gamma(X, \xi; \tau_1,...,\tau_n) := X + \int \frac{\delta}{\xi - \delta c \int \frac{\tau}{\tau X + 1}dH(\tau)}dD(\delta), \\
	&\hat X_i^j(\xi; \tau_1,...,\tau_n) := \check X(\xi).
\end{aligned}
\end{equation} 
We than have:
\begin{equation}\label{}
\begin{aligned}
	\frac{\partial \check X_i^j}{\partial \tau_k} = \frac{\partial \hat X_i^j}{\partial \tau_k} + \frac{\partial \check X_i^j}{\partial \xi_i^j}\frac{\partial \xi_i^j}{\partial \tau_k},
\end{aligned}
\end{equation} 
where:
\begin{equation}\label{}
\begin{aligned}
	& \frac{\partial \hat X_i^j}{\partial \tau_k} = - \frac{\frac{\partial \hat \Gamma_i^j}{\partial \tau_k}}{\frac{\partial \hat \Gamma_i^j}{\partial X}}, \\
	& \frac{\partial \check X_i^j}{\partial \xi_i^j} = - \frac{\frac{\partial \hat \Gamma_i^j}{\partial \xi_i^j}}{\frac{\partial \hat \Gamma_i^j}{\partial X}},
\end{aligned}
\end{equation} 
and:
\begin{equation}\label{}
\begin{aligned}
	& \frac{\partial \hat \Gamma_i^j}{\partial \tau_k} = \int \frac{\delta}{\left(\xi_i^j - \delta c \int \frac{\tau}{\tau \check X_i^j + 1}d\hat H(\tau)\right)^2} \times \frac{\delta c w_k}{(\tau_k \check X_i^j +1)^2}dD(\delta),\\
	& \frac{\partial \hat \Gamma_i^j}{\partial X}  = 1 - \int \frac{\delta}{\left(\xi_i^j - \delta c \int \frac{\tau}{\tau \check X_i^j + 1}d\hat H(\tau)\right)^2} \times \int\frac{\delta c\tau^2}{(\tau \check X_i^j +1)^2}d\hat H(\tau)dD(\delta),\\
	& \frac{\partial \hat \Gamma_i^j}{\partial \xi_i^j} = -\int \frac{\delta}{\left(\xi_i^j - \delta c \int \frac{\tau}{\tau \check X_i^j + 1}d\hat H(\tau)\right)^2}dD(\delta).
\end{aligned}
\end{equation} 

Finally, we are only interested in the asymptotic density, which is null on border of the spectrum, so for $j=0$ and $j=\omega_i+1$, just define $\check X_i^j = 0$ and $\frac{\partial \check X_i^j}{\partial \tau_k} = 0$.

\subsubsection{Step (4)}
In this step, we compute the asymptotic spectral density $f_i^j := \tilde F'(\xi_i^j)$. 

Let $i \in \llbracket 1, \nu \rrbracket$. For $j=0$ and $j=\omega_i+1$, we have $f_i^j = 0$ and $\frac{\partial f_i^j}{\partial \tau_k} = 0$.

For $j \in \llbracket 1, \omega_i \rrbracket$, we have:
\begin{equation}\label{}
\begin{aligned}
	& f_i^j = \frac{1}{\pi} \Im\left[-\frac{1}{\xi_i^j} \int \frac{1}{\tau \check X_i^j + 1}d\hat H(\tau)\right].
\end{aligned}
\end{equation} 
For the derivatives, we define:
\begin{equation}\label{}
\begin{aligned}
	& \tilde f(X, \xi; \tau_1,...,\tau_n) = -\frac{1}{\xi} \int \frac{1}{\tau \check X + 1}d\hat H(\tau).
\end{aligned}
\end{equation} 
We have:
\begin{equation}\label{}
\begin{aligned}
	& \frac{\partial f_i^j}{\partial \tau_k} =  \frac{1}{\pi} \Im\left[\frac{\partial \tilde f_i^j}{\partial X}\frac{\partial \check X_i^j}{\partial \tau_k} + \frac{\partial \tilde f_i^j}{\partial \xi}\frac{\partial \xi_i^j}{\partial \tau_k} + \frac{\partial \tilde f_i^j}{\partial \tau_k}\right],
\end{aligned}
\end{equation}
with:
\begin{equation}\label{}
\begin{aligned}
	& \frac{\partial \tilde f_i^j}{\partial X} = \frac{1}{\xi_i^j} \int \frac{\tau}{\left(\tau \check X_i^j + 1\right)^2}d\hat H(\tau),\\
	& \frac{\partial \tilde f_i^j}{\partial \xi} = \frac{1}{{\xi_i^j}^2} \int \frac{1}{\tau \check X_i^j + 1}d\hat H(\tau),\\
	& \frac{\partial \tilde f_i^j}{\partial \tau_k} = \frac{1}{\xi_i^j} \frac{w_k \check X_i^j}{\left(\tau_k \check X_i^j + 1\right)^2}.
\end{aligned}
\end{equation}

We can integrate $\tilde F'$ in order to retrieve $\tilde F$. We have the following continuity conditions:
\begin{equation}\label{}
\begin{aligned}
	& \tilde F(0) = \max(0, 1 - c^{-1}), \\
	& \tilde F_i^{\omega_i+1} :=  \tilde F(\xi_i^{\omega_i+1}) = \tilde F(0) + \left(1-\tilde F(0)\right)\sum_{j=1}^i \tilde \omega_j, \\
	& \tilde F_i^{\omega_i+1}  = \tilde F_i^0.
\end{aligned}
\end{equation}
Furthermore, we will integrate $f_i^j$ through trapezoidal integration, that will give us $G_i^j$. Then, we will rescale $G_i^j$ for each $i \in \llbracket 1, \nu \rrbracket$ so that it respects the border conditions of step (2), $\tilde F_i^{\omega_i+1} - \tilde F_i^0 = \tilde \omega_i$: this will be $\tilde F_i^j$.

For $j \in \llbracket 0, \omega_i+1 \rrbracket$, we have then:
\begin{equation}\label{}
\begin{aligned}
	& G_i^j = \tilde F_i^0 + \frac{1}{2}\sum_{\ell = 1}^j (\xi_i^\ell - \xi_i^{\ell-1})(f_i^\ell + f_i^{\ell-1}), \\
	& \tilde F_i^j = \tilde F_i^0 + \left(1-\tilde F(0)\right)\frac{\tilde \omega_i(G_i^j - \tilde F_i^0)}{G_i^{\omega_i+1} - \tilde F_i^0}.
\end{aligned}
\end{equation}

The derivatives are:
\begin{equation}\label{}
\begin{aligned}
	& \frac{\partial G_i^j}{\partial \tau_k} = \frac{1}{2}\sum_{\ell = 1}^j \left(\frac{\partial \xi_i^\ell}{\partial \tau_k} - \frac{\partial \xi_i^{\ell-1}}{\partial \tau_k}\right)(f_i^\ell + f_i^{\ell-1}) + \frac{1}{2}\sum_{\ell = 1}^j (\xi_i^\ell - \xi_i^{\ell-1})\left(\frac{\partial f_i^\ell}{\partial \tau_k} + \frac{\partial f_i^{\ell-1}}{\partial \tau_k}\right), \\
	& \frac{\partial \tilde F_i^j}{\partial \tau_k} = \left(1-\tilde F(0)\right)\frac{\tilde \omega_i}{G_i^{\omega_i+1}  - \tilde F_i^0}\frac{\partial G_i^j}{\partial \tau_k} -\left(1-\tilde F(0)\right)\frac{\tilde \omega_i (G_i^j - \tilde F_i^0)}{(G_i^{\omega_i+1} - \tilde F_i^0)^2}\frac{\partial G_i^{\omega_i+1}}{\partial \tau_k}.
\end{aligned}
\end{equation}

\subsubsection{Step (5)}
This last step consists in computing the loss between $\tilde F$ and the observed $F$. For this step, we compared three different loss, and used auto-differentiation in PyTorch to automatically compute the derivatives with respect to the inputs. This functionality gives the practitioner some flexibility on the choice of the suitable loss, as implementing the derivatives is not necessary. We considered:
\begin{itemize}
	\item the Wasserstein loss $\ell_{\text{Wass}}(\tilde F, F) = \lVert \tilde F - F \rVert_{\mathcal{W},2}^2$, as we are in dimension one, this is just the $2$-norm of the difference of the c.d.f.s,
	\item the quantile loss with linear interpolation $\ell_{\text{rect}}(\tilde F, F) = \sum_{i=1}^n (\tilde \lambda_i - \lambda_i)^2$, where $\tilde \lambda_i$ is $\frac{i}{n}^\text{th}$ quantile of $\tilde F$, approximated with linear interpolation on the grid $(\xi_i)_i$, and the sample eigenvalue $\lambda_i$ is by definition the $\frac{i}{n}^\text{th}$ quantile of the observed $F$,
	\item the quantile loss with trapezoidal integration $\ell_{\text{trap}}(\tilde F, F) = \sum_{i=1}^n (\tilde \lambda_i - \lambda_i)^2$ as proposed in \cite{Ledoit2016}, where $\tilde \lambda_i$ is this time approximated by $Z\left(\frac{i}{n}\right) - Z\left(\frac{i-1}{n}\right)$ with $Z: \kappa \mapsto \int_0^\kappa \tilde F^{-1}(u)du$, and $Z$ is itself approximated through trapezoidal integration.
\end{itemize}
We use two additional loss terms to make the optimization more robust:
\begin{itemize}
	\item the spectral mean norm: $\rho_1 \left(\sum_{i=1}^n \tilde \lambda_i - \hat \tau_i\right)^2$, motivated by the fact that $\int \lambda dF(\lambda) = \int \tau dH(\tau)$,
	\item a weight scheme for $\ell_{\text{rect}}$ and $\ell_{\text{trap}}$ of the form:
		\begin{equation*}\label{}
		\begin{aligned}
			&\ell_{\rho_2}(\tilde F, F) = \rho_2 \sum_{i=1}^n \frac{q_i}{\sum_{j=1}^n q_j/n} (\tilde \lambda_i - \lambda_i)^2,\text{ with } q_i = \cos\left(\frac{\pi i}{n}\right)^2.
		\end{aligned}
		\end{equation*} 
		The motivation is to increase the focus on extreme eigenvalues, which are of particular importance in linear applications of the covariance and precision matrix.
\end{itemize}

By default, our implementation uses the loss $\ell_{\text{trap}}$ with additional regularization parameters $\rho_1 = 0.1$ and $\rho_2 = 1$. The loss $\ell_{\text{trap}}$ is what we found to be the most numerically stable and suitable for the minimization, and the regularization helps at the extrema of the support.

The full \textit{WeSpeR} algorithm is then:
\begin{itemize}
	\item[1-] As input, we take the observed sample spectrum distribution $F_n = \frac{1}{n} \sum_{i=1}^n 1_{[\lambda_i,\infty[}$ and the weight matrix $W$.
	\item[2-] Find the estimated population spectrum $\hat H(\tau) = \frac{1}{n}\sum_{i=1}^n 1_{[\tau_i,+\infty[}$ where $\tau$ solves:
		\begin{equation}\label{wesper}
		\begin{aligned}
			\min_{\tau \in \R^n} \ell\left(\tilde F(c, \hat H) - F_n\right)
		\end{aligned}
		\end{equation} 
	using Adam optimizer.
	\item[3-] Compute $\check X(\lambda_i)$ with Proposition 4 \cite{Oriol2025b}, solving if $\lambda_i \in S_F$: 
		\begin{equation}
		\begin{aligned}
			\check X(\lambda_i) = \argmin_{X \in \C_+} \left| X + \int \frac{\delta}{\lambda_i - \delta c \int \frac{\tau}{\tau X + 1}d\hat H(\tau)}dD(\delta) \right|^2.
		\end{aligned}
		\end{equation} 
	\item[4-] Compute $h(\lambda_i)$ using Equation \ref{nl_formulas} with $\hat H$ and $\check X(\lambda_i)$.
\end{itemize}
The complete implementation (along with the (LD) and (HD) approaches) is provided in this \href{https://www.github.com/nlcvbo/WeSpeR}{Github repository}.

\begin{figure}[]
	\centering
	\includegraphics[width=0.49\linewidth]{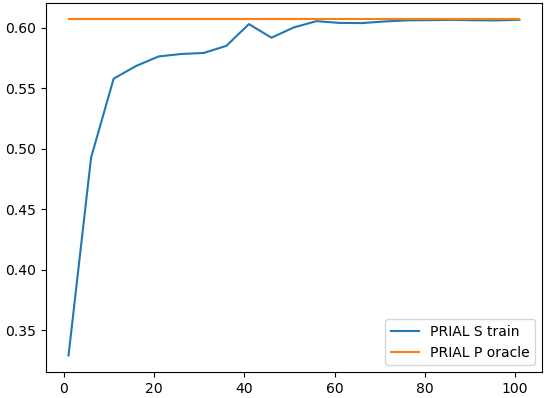}
	\includegraphics[width=0.49\linewidth]{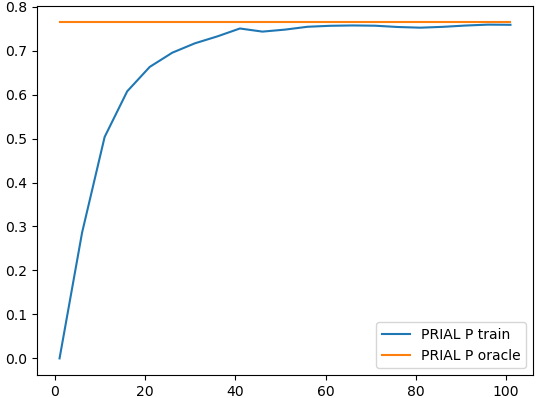} \\ 
	\includegraphics[width=0.49\linewidth]{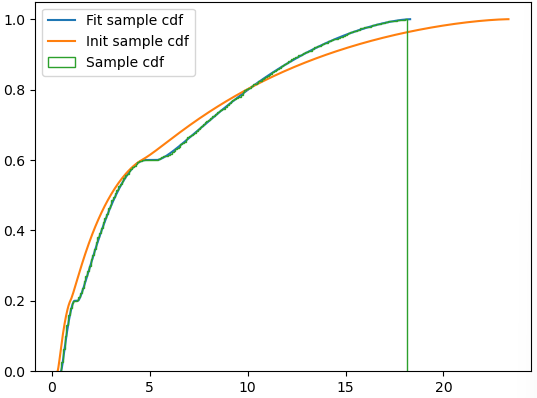}
	\includegraphics[width=0.49\linewidth]{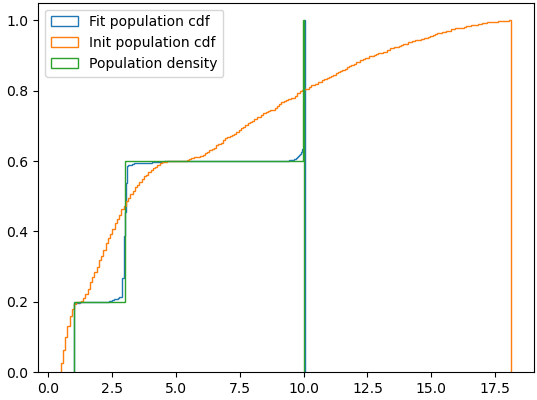}
	\includegraphics[width=0.7\linewidth]{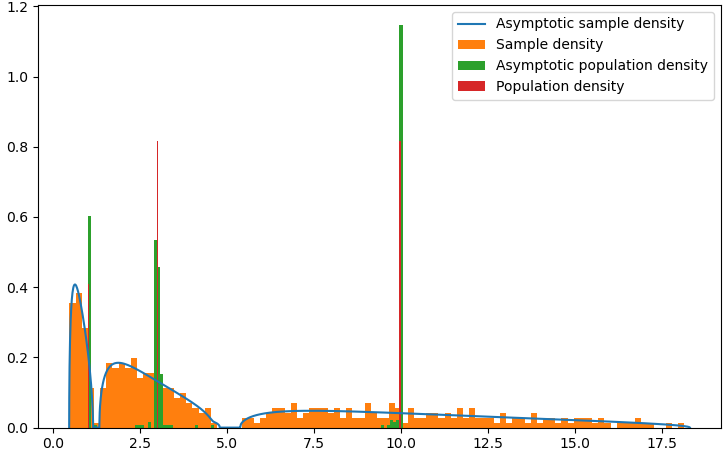}
	\caption[Result of \textit{WeSpeR} (MD) algorithm for $D$ exponentially-weighted with $\alpha=5.$]{Result of \textit{WeSpeR} algorithm for $D$ exponentially-weighted with $\alpha=5.$, $c=0.1$ and $n=400$ for readability. On top is displayed the PRIAL for covariance (left) and precision matrix (right) compared to the oracle epoch by epoch of the minimization step. On middle-left is displayed the initial and final sample c.d.f. $\tilde F$ after fitting along with the observed one $F_n$. On middle-right, we plotted the initial and final population c.d.f. $\hat H$ after fitting along with the population one $H$. On the bottom plot is shown the asymptotic sample density $\tilde F'$ estimated by \textit{WeSpeR} (blue), the histograms of the observed sample eigenvalues (orange), along with the true population spectrum histogram $H$ (red) and the estimated population spectrum histogram $\hat H$ (green).}
	\label{fig:wesper}
\end{figure}

An example of use of the algorithm in displayed in Figure \ref{fig:wesper}. We recall that the PRIAL of an estimator $\hat S$ for covariance estimation is defined as, with $\Sigma$ the population covariance and $S$ the sample covariance:
\begin{equation}
\begin{aligned}
	\text{PRIAL}_\text{cov} = 1 - \frac{\E\left[\lVert \hat S - \Sigma \rVert_F^2\right]}{\E\left[\lVert S - \Sigma \rVert_F^2\right]}.
\end{aligned}
\end{equation}
Similarly, for a precision matrix estimator $\hat P$, we have, with $S^\dagger$ the Moore-Penrose inverse:
\begin{equation}
\begin{aligned}
	\text{PRIAL}_\text{prec} = 1 - \frac{\E\left[\lVert \hat P - \Sigma^{-1} \rVert_F^2\right]}{\E\left[\lVert S^\dagger - \Sigma^{-1} \rVert_F^2\right]}.
\end{aligned}
\end{equation}

\begin{figure}[]
	\centering
	\includegraphics[width=0.49\linewidth]{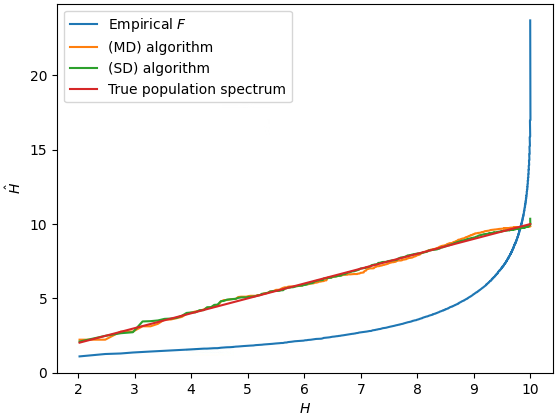} 
	\includegraphics[width=0.49\linewidth]{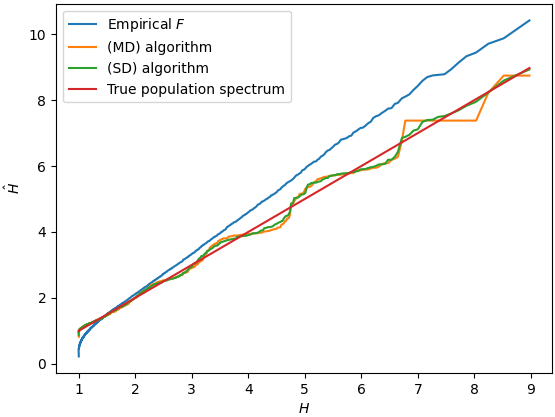} \\
	\includegraphics[width=0.49\linewidth]{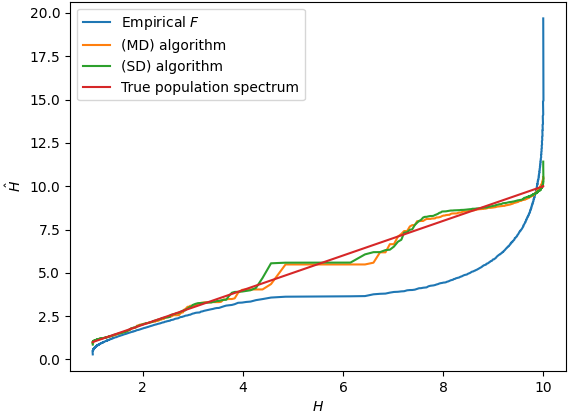} 
	\includegraphics[width=0.49\linewidth]{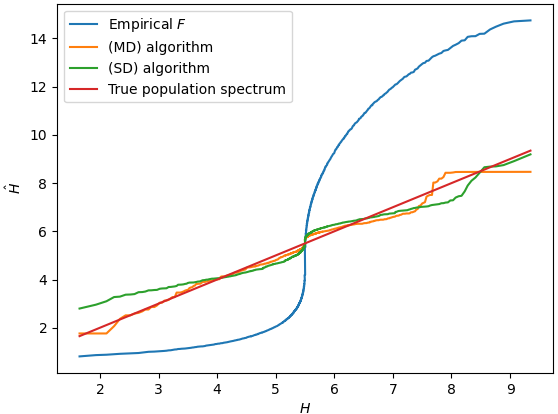} \\
	\caption[-Q plots between the true population spectrum $H_k$ and the estimated one $\hat H_k$ using (LD) and (MD) procedures.]{Q-Q plots between the true population spectrum $H_k$ and the estimated one $\hat H_k$ using (LD) and (MD) procedures for $100$ epochs, top-left for $H_1$, top-right for $H_2$, bottom-left for $H_3$ and bottom-right for $H_4$.}
	\label{fig:qq}
\end{figure}

We conducted additional experiments creating Q-Q plots, with the same setting as in \cite{Ledoit2016}, Section 10.1: four different population spectrum $H_k = 1 + (\eta - 1)G_k$ for $k \in \{1, 2, 3, 4\}$, with $\eta = 10$, $c=1/3$ and normal noise. For $x \in [0,1]$, $G_k$ are defined as follows:
\begin{equation}
\begin{aligned}
	&G_1(x) = 1 - (1-x^3)^{1/3}, \\
	&G_2(x) = [1-(1-x)^3]^{1/3}, \\
	&G_3(x) = \begin{cases}
					\frac{1}{2} \left[1 - (1-2x)^3\right]^{1/3} & \text{ if } x \in [0,1/2],\\
					1 - \frac{1}{2} \left[1 - (2x-1)^3\right]^{1/3} & \text{ if } x \in [1/2, 1],
				\end{cases} \\
	&G_4(x) = \begin{cases}
					\frac{1 - \left[1 - (2x)^3\right]^{1/3}}{2}  & \text{ if } x \in [0,1/2],\\
					\frac{1 + \left[1 - (2-2x)^3\right]^{1/3}}{2} & \text{ if } x \in [1/2, 1],
				\end{cases}
\end{aligned}
\end{equation}
These distributions include both null and infinite densities, and are good example to see the behavior of the algorithms in hard conditions. We introduce an additional weight distribution $D = \mathcal{U}([1/2,3/2])$ to this setting and retrieve the population spectrum from the weighted sample covariance in dimension $n = 1000$. We provide the Q-Q plots between the population spectrum $H_k$ and the estimated one $\hat H_k$ for each $k \in \llbracket 1, 4 \rrbracket$, plotting also the sample spectrum $F$ observed on the weighted sample covariance. The results are shown in Figure \ref{fig:qq} and highlight the fitting power of the (SD) and (MD) procedures. Constant regions mainly appear in low-density intervals where no population eigenvalues where sampled with $n=1000$ there.

Regarding the complexity, in step (2), each minimization step now only requires $O(nKM + n \kappa)$ time complexity, where $M$ is the number of connected components of $S_D$ and $K$ is the number of connected components of $\hat H$. Using the heuristic of Remark \ref{KD}, particularly useful when $D$ or $H$ is a mixture of Dirac with no large gaps, the step complexity falls to $O(n \kappa)$.

This complexity makes it scalable to large dimension, however for low dimension $n < 1000$ its performance is brought back by a large coefficient factor compared to the (LD) approach. This is mainly due to the fact that finding $S_{\tilde F}$ and $\tilde F$ is done using optimization routines or nested optimization routines.

\begin{figure}[]
	\centering
	\includegraphics[width=0.6\linewidth]{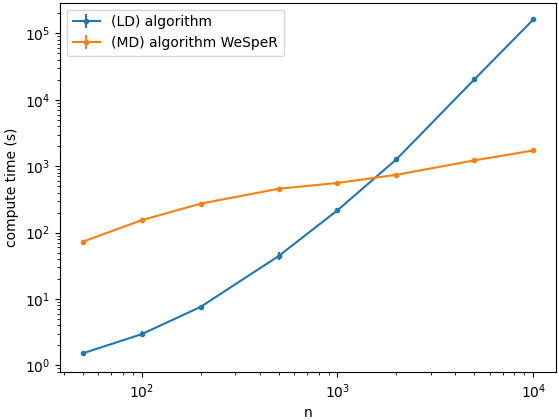} \\
	\includegraphics[width=0.49\linewidth]{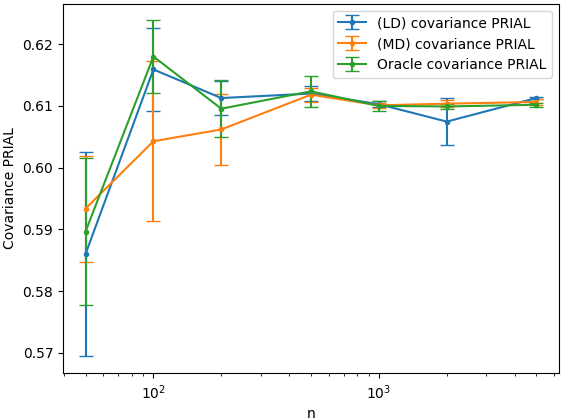}
	\includegraphics[width=0.49\linewidth]{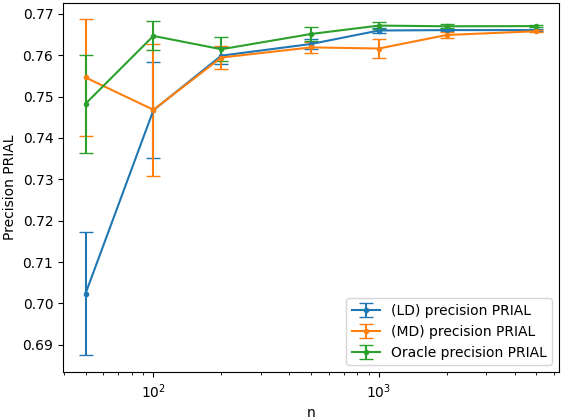}
	\caption[Compute time and PRIAL comparison for covariance and precision shrinkage between the (LD) and (MD) approaches, and the oracle.]{Compute time (top) and PRIAL comparison (bottom) for covariance (left) and precision (right) shrinkage between the (LD) and (MD) approaches, and the oracle. We fixed the number of epochs to $N_e = 100$, the weight distribution to a EWMA with parameter $\alpha = 5$, $c=0.1$, $\kappa = \min(n,400)$ and $H = \frac{1}{5} \mathbf{1}_{[1.,+\infty[} + \frac{2}{5} \mathbf{1}_{[3.,+\infty[} + \frac{2}{5} \mathbf{1}_{[3.,+\infty[}$. The experiment was run $11$ times, the error bars are $\pm \sigma$.}
	\label{fig:compute}
\end{figure}

In order to give order of magnitude of compute time and a performance comparison between the (LD) and (MD) approaches, we run the algorithms on a computer with $64$ GB of RAM and Intel Xeon(R) E-2286M CPU for different dimension $n$. The results and experiment details are given in Figure \ref{fig:compute}. The experiment shows us that while the compute time cross around $n = 2000$, the PRIAL stays similar and more importantly close to the oracle one, obtained using the true population spectrum $H$.

\section{\textit{WeSpeR} in the high dimensional (HD) case}
In the previous (MD) approach, using the \textit{WeSpeR} algorithm to compute the non-linear shrinkage covariance matrix estimator $S^*_n = U_n h(\Lambda) U_n^*$, we still use as input the observed eigenvalues $\lambda_i$ of $S_n$ and build $S_n^*$ with the diagonalization of $S_n$. This operation costs $O(n^3)$ once, at the beginning, and is what will make the approach unusable for very large dimension $n > 20000$ roughly. For reference, diagonalizing a $20000 \times 20000$ on a computer with $64$ GB of RAM and Intel Xeon(R) E-2286M CPU takes around $1000$ seconds, and computing $t(\lambda)$ or $h(\lambda)$ for all the eigenvalues takes also around $520$ seconds. Hopefully, we can overcome this limit with Lanczos algorithm \cite{Lanczos1950}.

For the high range of dimension $n > 20000$, we focus on linear applications of the covariance or precision matrices. Formally, we consider the problem of computing $f(S_n^*) v$ or $f(P_n^*) v$ for some vector $v \in \C^n$ and real function $f: \R \rightarrow \R$. By convention, we use the notation $f(M) := U f(\Lambda) U^*$ for $M = U \Lambda U^*$ in diagonalized form. This setting include linear regression or Markowitz optimization for example.

This comes to compute $(f \circ h)(S_n) v$ for covariance shrinkage, and $(f \circ t)(S_n) v$ for precision shrinkage. To compute it efficiently, we are using Stochastic Lanczos Quadrature (SLQ) and Lanczos algorithm \cite{Lanczos1950}. SLQ is analyzed for spectrum approximation in \cite{Chen2021}, where they give an upper bound of the number of vectors $n_v > 4(n+2)^{-1}t^{-2}\log(2n\eta^{-1})$ needed in SLQ for a given level of approximation $t$ of $F_n$ with given probability $\eta$. The running time is $O(n^2 t^{-1}\log(t^{-2}\eta^{-1}))$ in our case where $S_n$ is supposed to be dense \cite{Chen2010}.

The variant of \textit{WeSpeR} in the high dimension case is the following:
\begin{itemize}
	\item[1-] Estimate $F_n$ from $S_n$ using SLQ with $n_v$ vectors.
	\item[2-] Run \textit{WeSpeR} with $\kappa$ discretization points to estimate $\hat H$.
	\item[3-] Use Lanczos algorithm to compute $(f \circ h)(S_n) v$ (or $(f \circ t)(S_n) v$), where $h$ (or $t$) can be computed with $\hat H$, $c$ and $W$ on each eigenvalue in the Lanczos approximation.
\end{itemize}

The overall time complexity of the algorithm, with $N_{ep}$ epochs in the optimization step of \textit{WeSpeR} to estimate $\hat H$, is $O(n^2 t^{-1}\log(t^{-2}\eta^{-1}) + n\kappa N_{ep})$ using Remark \ref{KD}. In practice, the optimization part $O(n\kappa N_{ep})$ of \textit{WeSpeR} can be long to compute due to the optimization routines. In order to reduce this cost, this is always possible to approximate $H$ by a mixture of $\tilde n < n$ Dirac instead of $n$. For a cost in performance, the optimization complexity is reduced to $O(\tilde n \kappa N_{ep})$. We used $\tilde n$ step in Lanczos algorithm too, leading to a complexity of $O( \tilde n n^2 + \tilde n \kappa N_{ep})$.

This (HD) variant of \textit{WeSpeR} can handle large $n > 20000$ efficiently. Remark that regarding the memory complexity of the algorithm, the \textit{WeSpeR} optimization is in $O(n)$ and SLQ and Lanczos algorithm only require computing $S_n u$ for some vector $u \in \C^n$, which is led by the storage of $S_n$. In practice, apart from storing $S_n$, which takes $8n^2$ bytes for float64, the algorithm runs with $O(n)$ memory complexity.

In order to demonstrate the compute time of the (HD) approach on a standard use case, here is the experiment setting of a shrunk linear regression:
\begin{itemize}
	\item We sample $S$ of size $(n,n)$ with diagonal population covariance $\Sigma$ so that the population spectrum is $H = \frac{1}{5} \mathbf{1}_{[1.,+\infty[} + \frac{2}{5} \mathbf{1}_{[3.,+\infty[} + \frac{2}{5} \mathbf{1}_{[3.,+\infty[}$, $c=0.5$ and EWMA weight distribution with $\alpha = 0.1$.
	\item We sample a vector $v \in \mathcal{S}^{n-1}(\R)$ uniformly on the sphere.
	\item We compute the oracle Ordinary Least Square estimator $w^* = \Sigma^{-1} v$, the oracle shrunk estimator $w = t(S) v$ using Lanczos algorithm and the true $t(\cdot)$ function, $\hat w = \hat t(S) v$ with the estimated $\hat t(\cdot)$ function from the (HD) procedure, and the naive weight estimator $w_{naive} = S^\dagger v$. 
	\item We compare the resulting weight vectors using MSE, and monitor the compute time and the performance with the PRIAL of an estimator $\hat w$ defined as $\text{PRIAL} = 1 - \frac{\E\left[\lVert \hat w - w^* \rVert^2\right]}{\E\left[\lVert w_{naive} - \Sigma^{-1} \rVert_F^2\right]}$.
\end{itemize}

\begin{figure}[]
	\centering
	\includegraphics[width=0.49\linewidth]{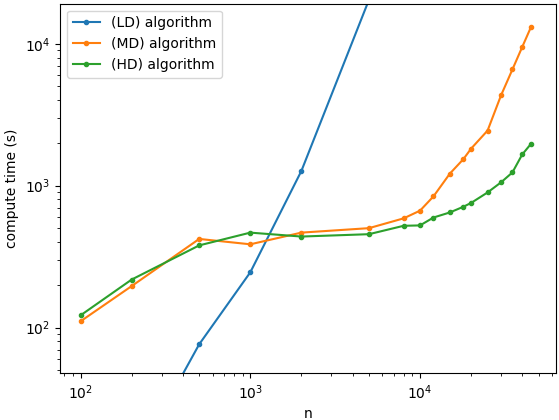}
	\includegraphics[width=0.49\linewidth]{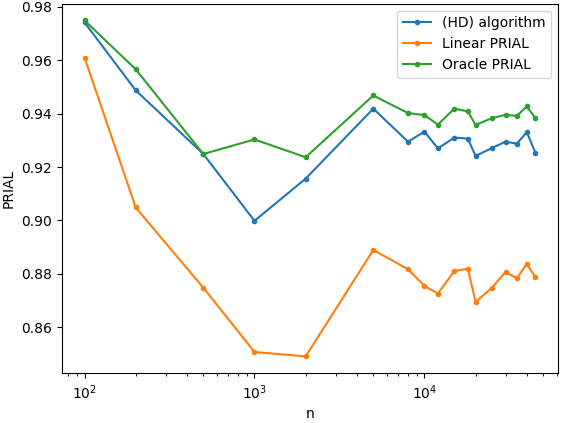}
	\caption[Compute time and PRIAL comparison for $w$ and $\hat w$.]{Compute time (top) and PRIAL comparison (bottom) for $w$ and $\hat w$. We fixed the number of epochs to $N_e = 100$, the weight distribution to a EWMA with parameter $\alpha = 0.1$, $c=0.5$, $\kappa = 400$, $n_v = 1$, $\tilde n=1000$, and $H = \frac{1}{5} \mathbf{1}_{[1.,+\infty[} + \frac{2}{5} \mathbf{1}_{[3.,+\infty[} + \frac{2}{5} \mathbf{1}_{[3.,+\infty[}$.}
	\label{fig:lanczos}
\end{figure}

The results are displayed in Figure \ref{fig:lanczos}. For comparison, we also show the compute time of (SD) and (MD) procedure described in previous sections. We clearly see the benefits of the two approximations proposed in the (HD) algorithm:
\begin{enumerate}
	\item we dodge a full diagonalization of $S_n$ in $O(n^3)$ using Lanczos algorithm in $O(\tilde n n^2)$,
	\item we reduce the cost of the optimization in \textit{WeSpeR} from $O(n\kappa N_{ep})$ to $O(\tilde n \kappa N_{ep})$.
\end{enumerate}
The advantage of this approach is to make it possible to trade compute time for precision and vice versa, while dodging a cubic time complexity arguably unusable in very large dimensions. Depending on the experimental setting, the available time and demand of performance, it is still possible to perform approximate non-linear shrinkage within reasonable compute time as long as the linear map $v \mapsto S_n v$ can be computed.

The performance in PRIAL are compared to the baseline, which is Ledoit-Wolf linear shrinkage \cite{Ledoit2004} adapted to weighted sample covariance. The associated PRIAL is denoted "Linear PRIAL".

Remark that in Figure \ref{fig:lanczos}, we used only $n_v=1$ vector for the SLQ. Indeed, in our experiment, using $n_v = 10$ would multiply the compute time by around $10$, while improving the PRIAL by less than $0.1\%$. This is a typical example of precision-compute time trade off where accepting a light diminishing of precision improves significantly the compute time here.

\section{Conclusion}
We have presented a comprehensive algorithmic framework for the computation of non-linear shrinkage formulas for weighted sample covariance matrices. In particular, we introduced \textit{WeSpeR}, a scalable and flexible algorithm that leverages theoretical properties of the asymptotic spectrum support. The algorithm is implemented using analytically computed derivatives as a PyTorch module. Empirical tests confirm the performance of the \textit{WeSpeR} algorithm in a wide range of dimensions. All three algorithms, (LD), (MD) and (HD) are implemented in Python in this \href{https://www.github.com/nlcvbo/WeSpeR}{Github repository}.

\newpage\section[Appendix - Experiments and implementation]{Appendix - Additional experiments, implementation details}

\subsection*{Support identification for $S_D$ convex}
\subsubsection{Examples of weight distribution: EWMA distribution}
\begin{figure}[ht]
\centering
\includegraphics[width=0.8\linewidth]{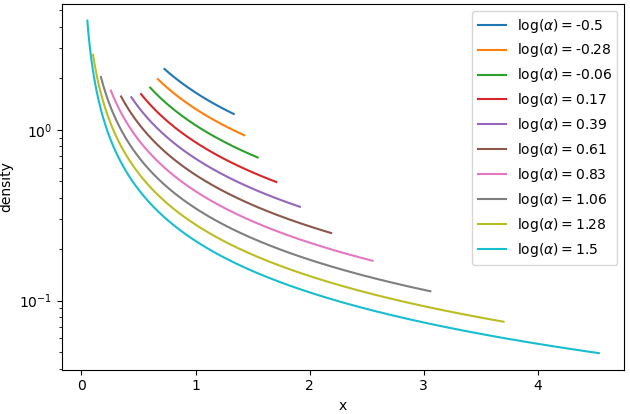}
\caption[Exponentially-weighted density $D_\alpha'$ for different $\alpha$.]{Exponentially-weighted density $D_\alpha'$ for different $\alpha$.}
\label{fig:exp}
\end{figure}

In time series analysis, the Exponentially-Weighted Moving Average - EWMA - is a widely used weighting scheme from neuroscience to finance as detailed in the introduction (see \cite{Pafka2004,Daly2010,Svensson2007,Tan2023}). The asymptotic c.d.f. $D_\alpha$ of the weights in a EWMA setting is defined in the following definition. The decay of the EWMA is controlled through the parameter $\alpha \in \R_+$: the larger $\alpha$, the steeper the decay. Some examples of densities for different values of $\alpha$ are shown in figure \ref{fig:exp}.

\begin{dft}[Exponentially-weighted (EWMA) distribution]\label{def:ewma}
	For $\alpha \in \R_+$, we consider the following weight c.d.f.:
	\begin{equation}\label{}
	\begin{aligned}
		D_\alpha: x \in [\beta e^{-\alpha}, \beta] \mapsto 1 + \frac{1}{\alpha}\log\left(\frac{x}{\beta}\right) \text{ with } \beta = \frac{\alpha}{1- e^{-\alpha}}.
	\end{aligned}
	\end{equation} 
	Then, we have the following closed-form formulas for $m_{LD}$ and $m_{LD}^{-1}$:
	\begin{equation}\label{}
	\begin{aligned}
		&m_{LD_\alpha}: x \in \R \backslash [\beta e^{-\alpha}, \beta] \mapsto \frac{1}{\alpha} \log \left(1 + \frac{\alpha}{\beta e^{-\alpha} - x} \right),\\
		&m_{LD_\alpha}^{-1}: y \in \R^*  \mapsto \beta e^{-\alpha} + \frac{\alpha}{1 - e^{\alpha y}}.
	\end{aligned}
	\end{equation} 
\end{dft}

\begin{figure}[]
\centering
\includegraphics[width=0.8\linewidth]{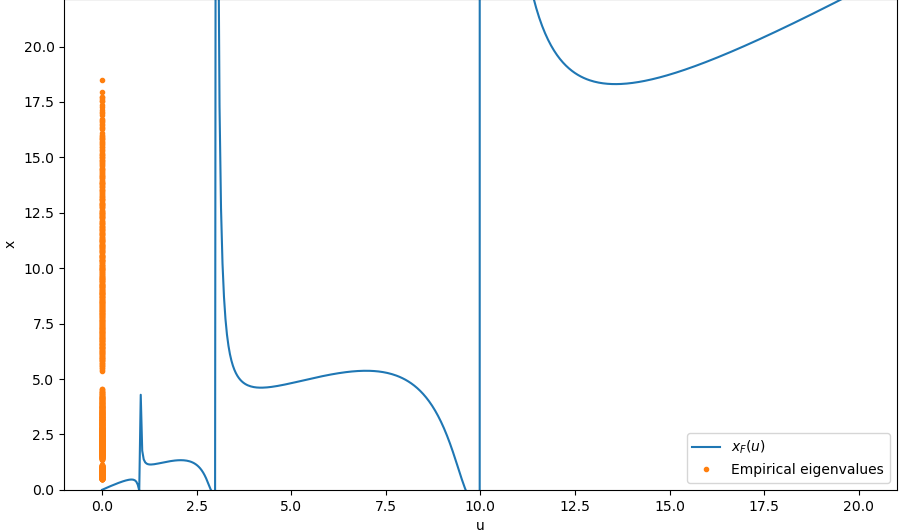}
\caption[Support identification procedure using $H$ mixture of $3$ Dirac, $D$ exponentially-weighted with $\alpha=5$, $c=0.1$.]{$x_F(-1/u)$ for $u \in S_H^c$, using $H$ being a mixture of 3 Dirac in $1$, $3$, and $10$ with respectively weights $0.2$, $0.4$ and $0.4$ as in \cite{Silverstein1995c}, $D_\alpha$ EWMA distribution with $\alpha = 5$, $c = 0.1$, and empirical eigenvalues sampled with $p=1000$. Horizontal lines are plotted at the zeros of $x_F'$, they represent the theoretical borders of $S_F$.}
\label{fig:sep_exp}
\end{figure}

\begin{figure}[]
\centering
\includegraphics[width=0.8\linewidth]{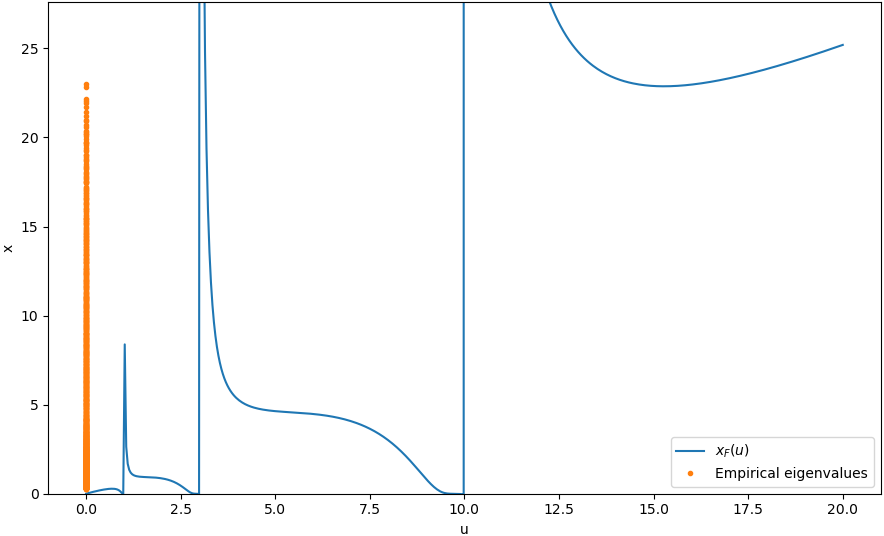}
\caption[Support identification procedure using $H$ mixture of $3$ Dirac, $D$ exponentially-weighted with $\alpha=10$, $c=0.1$.]{$x_F(-1/u)$ for $u \in S_H^c$, using $H$ being a mixture of 3 Dirac in $1$, $3$, and $10$ with respectively weights $0.2$, $0.4$ and $0.4$ as in \cite{Silverstein1995c}, $D_\alpha$ EWMA distribution with $\alpha = 10$, $c = 0.1$, and empirical eigenvalues sampled with $p=1000$. Horizontal lines are plotted at the zeros of $x_F'$, they represent the theoretical borders of $S_F$.}
\label{fig:nonsep_exp}
\end{figure}

Two examples of identification of the support are shown in figure \ref{fig:sep_exp} when there is spectral separation, and in figure \ref{fig:nonsep_exp} when there is not.

\subsubsection{Examples of weight distribution: uniform distribution}
\begin{ppt}[Uniform distribution]
For $\alpha \in [0, 2[$, we consider the following weight distribution:
\begin{equation}\label{}
\begin{aligned}
	D_\alpha: x \in \R \mapsto \frac{x - 1 + \alpha/2}{\alpha} 1_{[1-\alpha/2, 1 + \alpha/2]}(x) + 1_{[1+\alpha/2, +\infty[}(x).
\end{aligned}
\end{equation} 
Then, we have the following closed-form formulas for $m_{LD}$:
\begin{equation}\label{}
\begin{aligned}
	&m_{LD_\alpha}: x \in \R \backslash [1-\alpha/2, 1 + \alpha/2] \mapsto 1 + \frac{x}{\alpha} \log\left(1 + \frac{\alpha}{1 - \frac{\alpha}{2} - x} \right).
\end{aligned}
\end{equation} 
$m_{LD}^{-1}$ has no closed-form formulas, and it can be retrieved through numerical optimization.
\end{ppt}
\newpage

\begin{figure}[]
\centering
\includegraphics[width=0.75\linewidth]{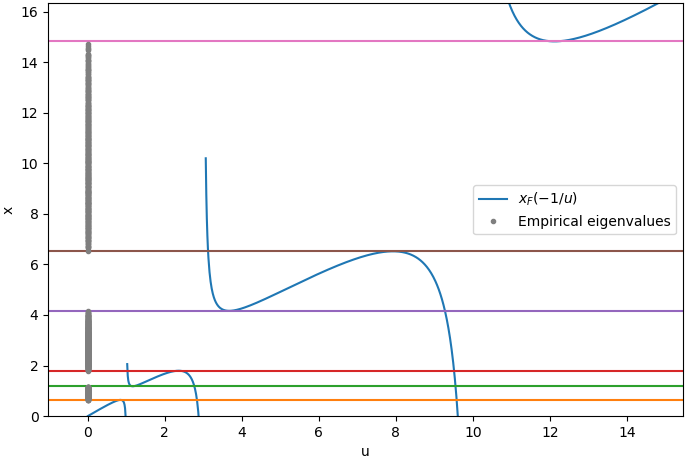}
\caption[Support identification procedure using $H$ mixture of $3$ Dirac, $D$ uniform with $\alpha=1$, $c=0.1$.]{$x_F(-1/u)$ for $u \in S_H^c$, using $H$ being a mixture of 3 Dirac in $1$, $3$, and $10$ with respectively weights $0.2$, $0.4$ and $0.4$ as in \cite{Silverstein1995c}, $D_\alpha$ uniform with $\alpha = 1$, $c = 0.1$, and empirical eigenvalues sampled with $p=1000$. Horizontal lines are plotted at the zeros of $x_F'$, they represent the theoretical borders of $S_F$.}
\label{fig:sep_unif}
\end{figure}

\begin{figure}[]
\centering
\includegraphics[width=0.8\linewidth]{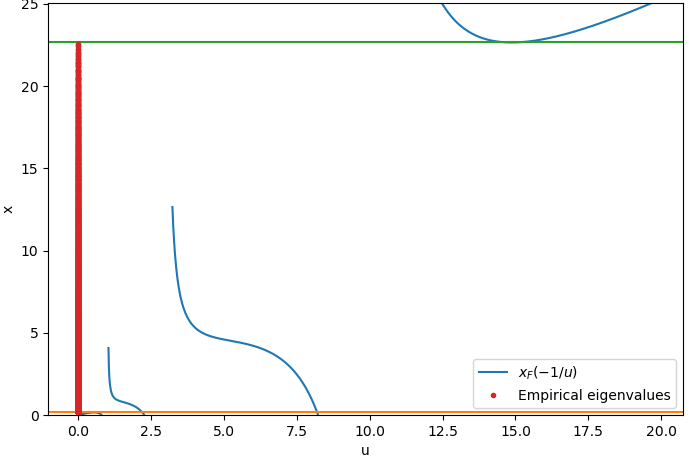}
\caption[Support identification procedure using $H$ mixture of $3$ Dirac, $D$ uniform with $\alpha=1$, $c=0.5$.]{$x_F(-1/u)$ for $u \in S_H^c$, using $H$ being a mixture of 3 Dirac in $1$, $3$, and $10$ with respectively weights $0.2$, $0.4$ and $0.4$ as in \cite{Silverstein1995c}, $D_\alpha$ uniform with $\alpha = 1$, $c = 0.5$, and empirical eigenvalues sampled with $p=5000$. Horizontal lines are plotted at the zeros of $x_F'$, they represent the theoretical borders of $S_F$.}
\label{fig:nonsep_unif}
\end{figure}

Two examples of identification of the support are shown in figure \ref{fig:sep_unif} when there is spectral separation, and in figure \ref{fig:nonsep_unif} when there is not.

\subsubsection{Numerical implementation of spectrum support identification}\label{implem}
In all the implementation, we assume we can compute $m_{LD}^{-1}$, $t: u \in S_H^c \mapsto c\int \frac{\tau}{\tau - u}dH(\tau)$, and the derivatives $m_{LD}'$, $m_{LD}''$, $t'$ and $t''$ at any point. Moreover, we assume we have access to a root-finding algorithm taking in input a real function $f$ and two real values $(x_0,x_1)$ such that $f(x_0)$ and $f(x_1)$ have different signs. The idea is inspired from the QuEST algorithm \cite{Ledoit2015}.

We implement the algorithm in the scenario where $S_H$ can be written as a finite disjoint union of intervals, \textit{i.e.} for some $K \in \N^*$, $S_H = \overset{K}{\underset{i=1}\cup} [\tau_i^l,\tau_i^r]$ where $\tau_1^l \leq \tau_1^r < ... < \tau_K^l \leq \tau_K^r$.

For convenience, we use the increasing change of variable $u \mapsto -1/X$, and consider the function $y_F: u \in S_H^c \mapsto x_F(-1/u)$. We have, for $u \in S_H^c$ with $t: u \in S_H^c \mapsto c\int \frac{\tau}{\tau - u}dH(\tau)$:
\begin{equation}\label{deriv}
	\small
\begin{aligned}
	y_F(u) = &-ut(u)m_{LD}^{-1}(t(u)),\\
	y_F'(u) = &-t(u)m_{LD}^{-1}(t(u)) -ut'(u)\left(\frac{t(u)}{m_{LD}'\left(m_{LD}^{-1}(t(u))\right)} + m_{LD}^{-1}(t(u)) \right),\\
	y_F''(u) = &-(2t'(u) + ut''(u)))m_{LD}^{-1}(t(u)) - \frac{2\left(t(u)+ut'(u)\right)t'(u) + ut(u)t''(u)}{m_{LD}'\left(m_{LD}^{-1}(t(u))\right)} \\
	&- \frac{ut(u)t'(u)^2m_{LD}''\left(m_{LD}^{-1}(t(u))\right)}{m_{LD}'\left(m_{LD}^{-1}(t(u))\right)^3}.
\end{aligned}
\end{equation} 

We are going to construct iteratively $S_F^c$. At step $0$, we consider $(S_F^c)_0 = \emptyset$. There are three different situations. 
\begin{itemize}
	\item Firstly, on the interval $]-\infty, \tau_1^l[$, we are looking for the unique zero of $y_F'$. As $y_F'(u) \underset{u \rightarrow -\infty}{\longrightarrow} \int \delta dD(\delta) > 0$ and $y_F'(u) \underset{u \rightarrow \tau_1^{l-}}{\longrightarrow} -\infty$, we can find easily with line search two points $(u_1,u_2) \in ]-\infty, \tau_1^l[^2$ such that $y_F'(u_1) > 0$ and $y_F'(u_2) < 0$. We can use the root-finding algorithm of $y_F'$ between $u_1$ and $u_2$, giving us the solution $u^*_l$. In conclusion of this part, we update $(S_F^c)_1 = (S_F^c)_0 \cup ]-\infty, y_F(u_l^*)[$.
	
	\item Similarly, on the interval $]\tau_K^r,+\infty[$, we are looking for the unique zero of $y_F'$. As $y_F'(u) \underset{u \rightarrow +\infty}{\longrightarrow} \int \delta dD(\delta) > 0$ and $y_F'(u) \underset{u \rightarrow \tau_K^{r+}}{\longrightarrow} -\infty$, we can apply the previous procedure in $]\tau_K^r,+\infty[$ and find $u^*_r$, root of $y_F'$ in $]\tau_K^r,+\infty[$. In conclusion of this part, we update $(S_F^c)_2 = (S_F^c)_1 \cup ]y_F(u_r^*),+\infty[$.
	
	\item For each $i \in \llbracket 1,K-1 \rrbracket$, we consider the interval $]\tau_i^r,\tau_{i+1}^l[$. This time, we have $y_F''(u) \underset{u \rightarrow \tau_i^{r+}}{\longrightarrow} +\infty$ and $y_F''(u) \underset{u \rightarrow \tau_{i+1}^{l-}}{\longrightarrow} -\infty$. Still through line-search, we can use the root-finding algorithm on $y_F''$. We expect $y_F''$ to have only one zero, denoted by $u_0 \in ]\tau_i^r,\tau_{i+1}^l[$ on this interval. 
		\begin{itemize}
			\item If $y_F'(u_0) \leq 0$, there is no spectral gap to be found on this interval.
			\item Otherwise, we are looking for two zeros of $y_F'$: one on $]\tau_i^r, u_0[$ and one on $]u_0,\tau_{i+1}^l[$. As $y_F'(u) \underset{u \rightarrow \tau_i^{r+}}{\longrightarrow} -\infty$ and $y_F'(u) \underset{u \rightarrow \tau_{i+1}^{l-}}{\longrightarrow} -\infty$, we can apply the line-search and use the root-finding algorithm on each interval, outputting two solutions: $u_l^* \in]\tau_i^r, u_0[ $ and $u_r^* \in ]u_0,\tau_{i+1}^l[$. We update our complementary support: $(S_F^c)_{i+2} = (S_F^c)_{i+1} \cup ]y_F(u_l^*),y_F(u_r^*)[$.
		\end{itemize}
\end{itemize}
In the end, we have $S_F^c = (S_F^c)_{K+1}$.

\subsection*{Support identification for $D$ mixture of Dirac}
We give examples of applications of the Theorem \ref{spectrum2}, when $S_D$ is a finite union of intervals.
\subsubsection{Examples of weight distribution: mixture of two Dirac}
We detail the simpler case of $S_D$ as union of several intervals: $S_D$ as union of two points, \textit{i.e.} $D$ being a mixture of two Dirac.
\begin{ppt}[Mixture of two Dirac]
For $\alpha \in ]0, 1[, w \in ]0,1[$, we consider the following weight distribution, with $\beta = \frac{\alpha w}{1-w}$:
\begin{equation}\label{}
\begin{aligned}
	D_{\alpha,w} = w1_{[1 - \alpha, \infty[} + (1-w)1_{[1 + \beta, \infty[}.
\end{aligned}
\end{equation} 
Then, we have the following closed-form formulas for $m_{LD}$ and $m_{LD}^{-1}$:
\begin{equation}\label{}
\begin{aligned}
	&m_{LD_{\alpha,w}}^{(1)}: x \in ]1 - \alpha, 1+\alpha[ \mapsto \frac{w (1-\alpha)}{1 - \alpha - x} + \frac{(1-w)(1+\beta)}{1 + \beta - x} ,\\
	&m_{LD_{\alpha,w}}^{(2)}: x \in \R \backslash [1 - \alpha, 1+\alpha] \mapsto \frac{w (1-\alpha)}{1 - \alpha - x} + \frac{(1-w)(1+\beta)}{1 + \beta - x} ,\\
	&\left(m_{LD_{\alpha,w}}^{(1)}\right)^{-1}: y \in \R^*  \mapsto \frac{-b(y)^2 + \sqrt{b(y)^2 -4a(y)c(y)}}{2a(y)},\\
	&\left(m_{LD_{\alpha,w}}^{(2)}\right)^{-1}: y \in \R^*  \mapsto \frac{-b(y)^2 - \sqrt{b(y)^2 -4a(y)c(y)}}{2a(y)},
\end{aligned}
\end{equation} 
with for $y \in \R^*$:
\begin{equation}\label{}
\begin{aligned}
	&a(y) = y(1-w), \\
	&b(y) = 1-w-\left(1-2w(1-\alpha) + 1-\alpha\right)y,  \\
	&c(y) = (1-\alpha)\left(1-w(1-\alpha)\right) - (1-\alpha)\left(1-w(a-\alpha)\right).
\end{aligned}
\end{equation} 
\end{ppt}

\begin{figure}[]
\centering
\includegraphics[width=0.9\linewidth]{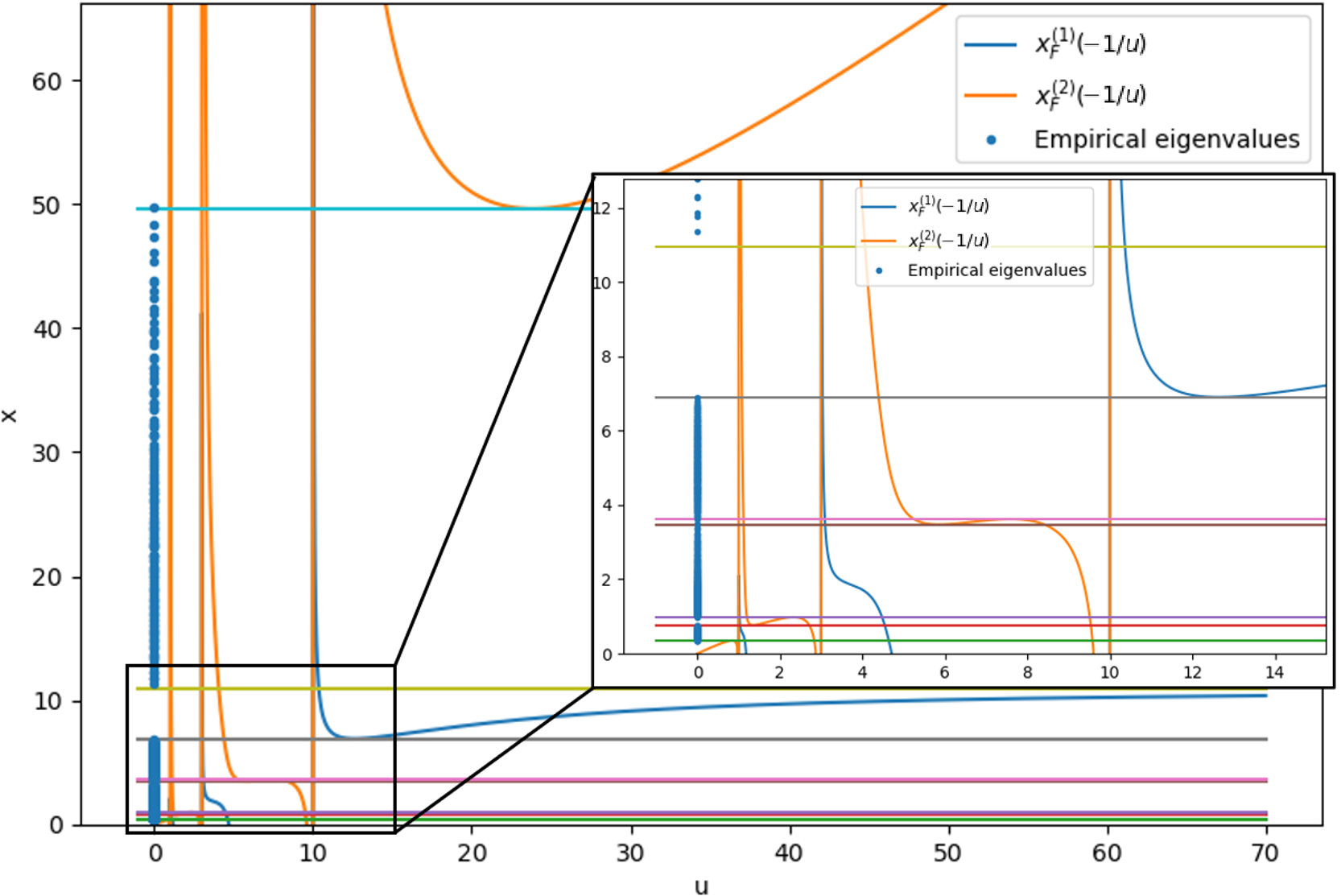}
\caption[Support identification procedure using $H$ mixture of $3$ Dirac, $D$ mixture of $2$ Dirac with $\alpha=0.5, w = 1 - \frac{1}{80}$, $c=0.1$.]{$x_F^{(k)}(-1/u)$ for $u \in S_H^c$, using $H$ being a mixture of 3 Dirac in $1$, $3$, and $10$ with respectively weights $0.2$, $0.4$ and $0.4$ as in \cite{Silverstein1995c}, $D_{\alpha,w}$ mixture of 2 Dirac with $\alpha = 0.5, w = 1 - \frac{1}{80}$, $c = 0.1$, and empirical eigenvalues sampled with $p=4000$. Horizontal lines are plotted at the zeros of $x_F^{(k)'}$, they represent the theoretical borders of $S_F$.}
\label{fig:sep_2Dirac2}
\end{figure}

An example of identification of the support is shown in figure \ref{fig:sep_2Dirac2} where a new spectral separation for the empirical spectrum is induced by the weight distribution - itself having large gaps. Thus, $S_F$ is made of $4$ intervals, while $S_H$ is only made of $3$. This is a new behavior due to the weight distribution $D$ that we do not observe in the classical setting with equal weights, where $S_F$ could only be made of at most $3$ intervals with this type of true spectrum $H$. 


\subsubsection{Examples of weight distribution: mixture of $N$ Dirac}
We detail the computation and the result for the $D$ mixture of $N$ Dirac. The efficient implementation of this problem is discussed after.
\begin{ppt}[Mixture of $N$ Dirac]
For $M \in \N^*, w \in ]0, 1]^M, \delta \in \left(\R_+^*\right)^M$ such that $\sum_{i=1}^M w_i = 1, \sum_{i=1}^M w_i \delta_i= 1$, we consider the following weight distribution:
\begin{equation}\label{}
\begin{aligned}
	D_\alpha = \sum_{i=1}^M w_i 1_{[\delta_i, \infty[}.
\end{aligned}
\end{equation} 
Then, we have the following closed-form formulas for $m_{LD}$:
\begin{equation}\label{}
\begin{aligned}
	& \forall k \in \llbracket 1, M-1 \rrbracket, m_{LD}^{(k)}: x \in ]\delta_2^{(k)},\delta_1^{(k+1)}[ \mapsto  \sum_{i=1}^M w_i \frac{\delta_i}{\delta_i - x}, \\
	& m_{LD}^{(M)}: x \in ]-\infty,\delta_1^{(1)}[ \cup ]\delta_2^{(M)},+\infty[ \mapsto  \sum_{i=1}^M w_i \frac{\delta_i}{\delta_i - x}.
\end{aligned}
\end{equation} 
$\left(m_{LD}^{(k)}\right)^{-1}$ has no general closed-form formulas when $M \geq 5$ due to Abel-Ruffini theorem, and it can be retrieved through numerical optimization.
\end{ppt}

\begin{figure}[]
\centering
\includegraphics[width=0.9\linewidth]{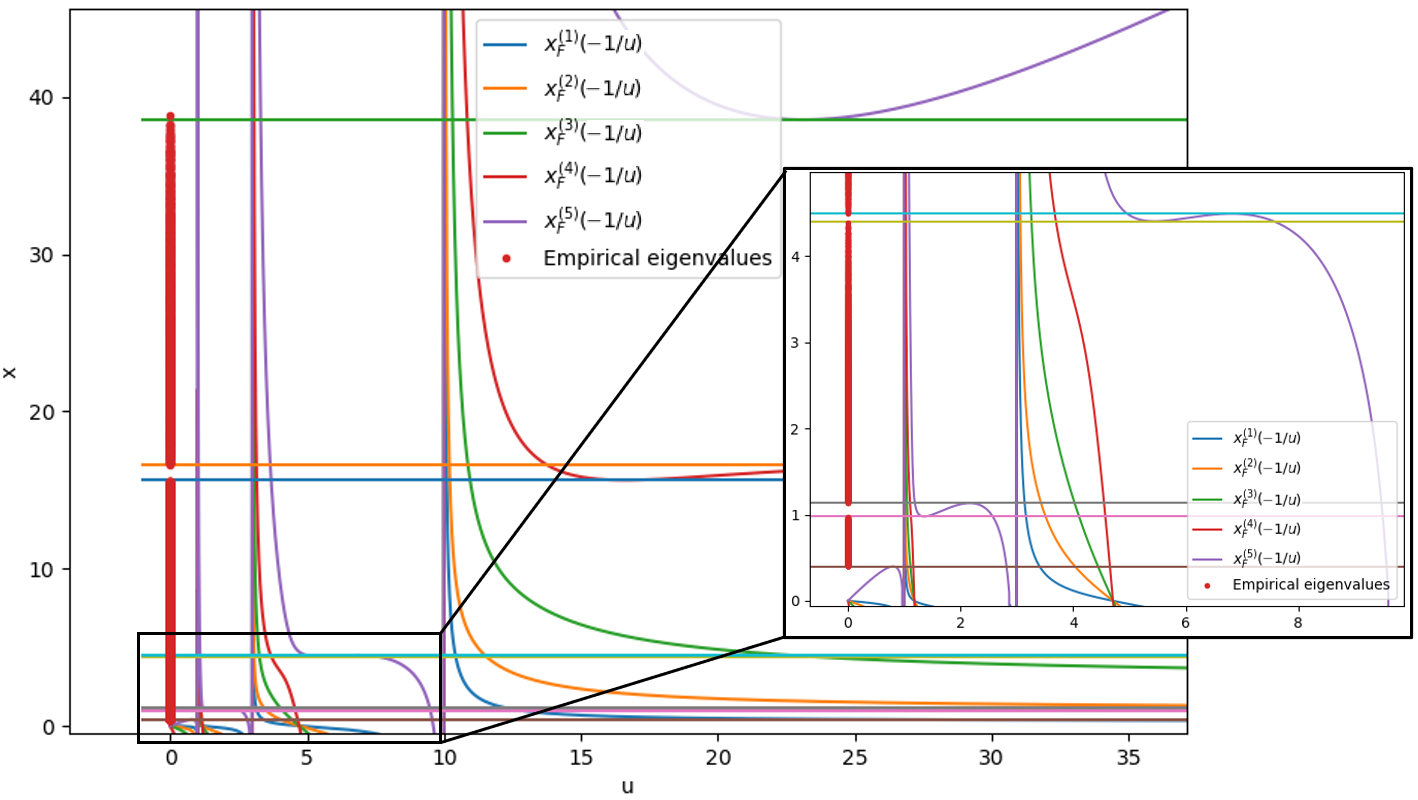}
\caption[Support identification procedure using $H$ mixture of $3$ Dirac, $D$ mixture of $5$ Dirac, $c=0.1$.]{$x_F^{(k)}(-1/u)$ for $u \in S_H^c$, using $H$ being a mixture of 3 Dirac in $1$, $3$, and $10$ with respectively weights $0.2$, $0.4$ and $0.4$ as in \cite{Silverstein1995c}, $D$ mixture of 5 Dirac with in $0.34, 0.67, 2.7, 6.74,$ and $34$, with respective weights $0.59, 0.30, 0.074, 0.03$, and $0.006$, $c = 0.1$, and empirical eigenvalues sampled with $p=5000$. Horizontal lines are plotted at the zeros of $x_F^{(k)'}$, they represent the theoretical borders of $S_F$.}
\label{fig:sep_NDirac}
\end{figure}

An example of identification of the support are shown in figure \ref{fig:sep_NDirac} when an important separation of the weight Dirac implies a new spectral separation for the empirical spectrum, in order to show the role of each $\left(m_{LD}^{(k)}\right)^{-1}$ in the determination of the support.

\subsubsection{Numerical implementation of spectrum support for $D$ finite mixture of Dirac}\label{implemNDirac}
We mostly use the same implementation scheme provided for the convex setting in Section \ref{implem}, applied to each $y_F^{(k)}: u \in S_F^c \mapsto x_F^{(k)}(-1/u)$ (extended in $0$ by $0$). 

Using the same notation, the only difference in the implementation is about the handling of $y_F^{(k)}$, $k \in \llbracket 1, M-1 \rrbracket$ on the outside interval $]\tau_K^r,+\infty[$. For numerical reasons, we apply back the change of variable $X = -1/u$, as it is easier to study $x_F^{(k)}$ on $]-1/\tau_K^r,0[$ than $y_F^{(k)}$ on the initial interval. On each of these intervals, we use the same idea we used on the bounded intervals: find the zero $X_0$ of $x_F^{(k)''}$, and if $x_F^{(k)'}(X_0)>0$, find the zero $X_l^*$ of $x_F^{(k)'}$ on $]-1/\tau_K^r,X_0[$ and the zero $X_r^*$ of $x_F^{(k)'}$ on $]X_0,0[$. Then add $]-1/x_F(X_l^*), -1/x_F(X_r^*)[$ to the current $(S_F^c)_i$.

The central numerical problem remaining is to compute the function $v^{(k)}: u \mapsto m_{LD}^{(k)-1}\left( c \int \frac{\tau}{\tau - u}dH(\tau)\right)$ for the desired $u \in S_H^c$. 

In this section, the term $t: u \in S_H^c \mapsto c\int \frac{\tau}{\tau - u}dH(\tau)$ is supposed to be easy to compute. In the case of $H$ being a finite mixture of Dirac, it is a rational function.

The more complex part resides in computing $v^{(k)}$. We suppose that $D$ is a finite mixture of Dirac, \textit{i.e.} there exists $M \in \N^*, (w_i)_{i=1}^M \in \R_+^*,  (\delta_i)_{i=1}^M \in \R_+^*$, such that $\sum_{i=1}^M w_i = 1$ and $D = \sum_{i=1}^M w_i 1_{[\delta_i,+\infty[}$.

Let $t \in \R$. Computing $v^{(k)}(t)$ for all $k \in \llbracket 1,M \rrbracket$ is equivalent to finding the $M$ distinct roots of the rational function $x \in S_D^c \mapsto m_{LD}(x) - t$. And this is equivalent to finding the $M$ distinct roots of the polynomial $P - tQ$ where:
\begin{equation}\label{}
\begin{aligned}
		&P(X) = \sum_{i=1}^M w_i \delta_i \prod_{j = 1,j \neq i}^M (\delta_j - X),\\
		&Q(X) = \prod_{i=1}^M (\delta_i - X).
\end{aligned}
\end{equation} 

We suggest to use \href{https://github.com/robol/MPSolve/tree/master}{MPSolve} (useable in Python, Matlab, C, Octave...) to find efficiently and simultaneously the $M$ roots of the resulting polynomial routinely for large $M$. Otherwise, for moderately large $M$, using the eigenvalues of the companion matrix of $P-tQ$ is possible and makes the implementation slightly easier.

Once we computed through this method all the functions $u \mapsto m_{LD}^{(k)-1}\left( c \int \frac{\tau}{\tau - u}dH(\tau)\right)$ on the desired grid of $u \in S_H^c$, we can easily deduce the $y_F^{(k)}$ and its derivative $y_F^{(k)'}$, $y_F^{(k)''}$ with the following formulas:
\begin{equation}\label{}
\small
\begin{aligned}
		y_F^{(k)}(u) = &-ut(u)v^{(k)}(t(u)),\\
		y_F^{(k)'}(u) = &-t(u)v^{(k)}(t(u)) -ut'(u)\left(\frac{t(u)}{m_{LD}'\left(v^{(k)}(t(u))\right)} + v^{(k)}(t(u)) \right),\\
		y_F^{(k)''}(u) = &-(2t'(u) + ut''(u)))v^{(k)}(t(u)) - \frac{2\left(t(u)+ut'(u)\right)t'(u) + ut(u)t''(u)}{m_{LD}'\left(v^{(k)}(t(u))\right)} \\ 
		&- \frac{ut(u)t'(u)^2m_{LD}''\left(v^{(k)}(t(u))\right)}{m_{LD}'\left(v^{(k)}(t(u))\right)^3}.
\end{aligned}
\end{equation} 

\begin{rmk}
	When only one root is necessary at a time, as it is in \textit{WeSpeR}, we recommend finding the roots of $m_{LD}$ directly without computing the polynomials, as it seems to be numerically more stable, as the coefficients of the polynomials can quickly overflow.
\end{rmk}

\subsection*{Additional experiments}
\begin{figure}[]
\centering
\includegraphics[width=0.57\linewidth]{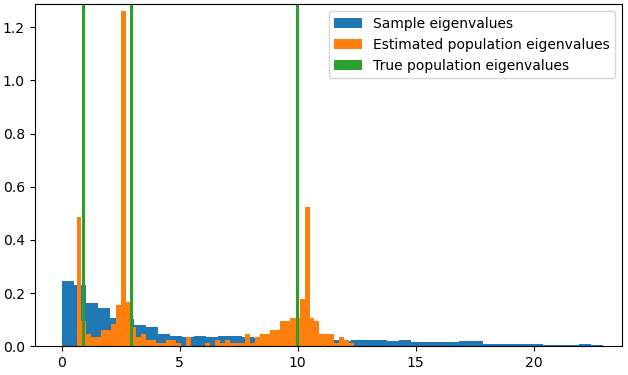}
\includegraphics[width=0.57\linewidth]{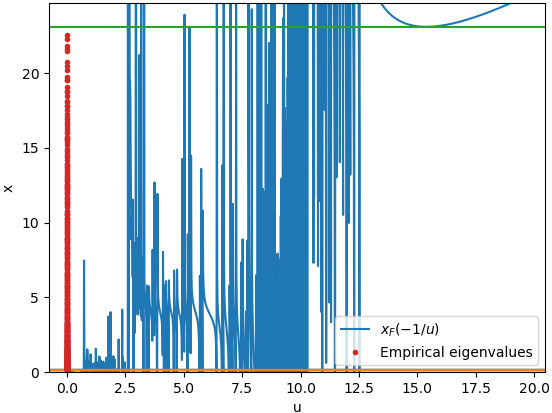}
\includegraphics[width=0.57\linewidth]{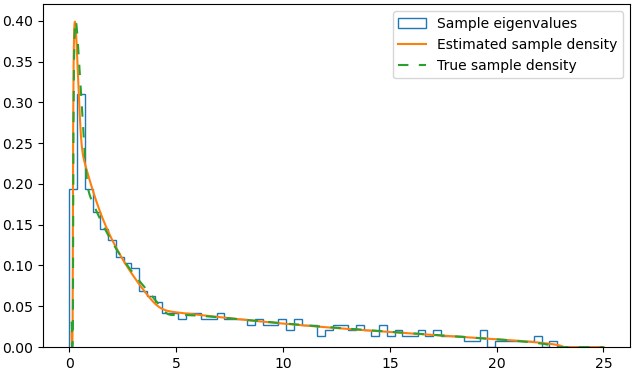}
\caption[Histograms, support identification and densities resulted from \textit{WeSpeR} (MD) procedure, Gaussian noise.]{(Top) Histograms of sample eigenvalues, estimated population eigenvalues $\hat H$ with \textit{WeSpeR} (MD), and true population eigenvalues $H$. (Middle) $u \mapsto x_F(-1/u)$ using the estimated $\hat H$ to detect the estimated support $S_F$. (Bottom) Estimated and true sample density computed on $S_F$ and sample eigenvalues' histogram. $H = \frac{1}{5} \mathbf{1}_{[1,\infty[} + \frac{2}{5} \mathbf{1}_{[3,\infty[} + \frac{2}{5} \mathbf{1}_{[10,\infty[}$, $D$ exponentially weighted with $\alpha = 1$, $c=0.5$, and $Z_{ij} \sim \mathcal{N}(0,1)$.}
\label{fig:wesper1}
\end{figure}

We explore the behavior of the algorithm with a different concentration ratio $c=0.5$ to see the effect of fewer samples on the estimation. Results are shown in Figure \ref{fig:wesper1}. Fewer samples compared to the dimension makes it hard to estimate accurately $H$. We clearly see that the estimated $H$ is more spread around the true Dirac than it were with $c=0.1$ in the experiment shown in the main corpus. 

We also experiment the effect of heavy tails in the estimation. We remark that even if the convergence of $F_n$ to $F$ is almost sure as long as we have bounded $2^{nd}$ moments, this convergence is slower as the tail is heavier. 

In this experiment, we fix the dimension $n = 400$, $c=0.1$, and we study the impact of $\nu > 2$ when the underlying noise $Z_{ij}$ follows a Student distribution $t_\nu(0,1)$. In this scenario, nothing particular happens while $\nu \geq 4$ roughly, $F_n$ is consistent, and the estimation remains barely affected. However, for very low $\nu$, around $2 < \nu < 4$, the sample spectrum $F_n$ tends to have some very high eigenvalues outside the asymptotic support $S_F$. Everything else fixed, decreasing $\nu$ increases the amount of outlier eigenvalues in $F_n$. Of course their frequency vanishes while $n \rightarrow +\infty$, but they exist in moderate dimension. 

In order to study the impact of these outliers in the estimation of $H$, we consider an extreme setting $\nu = 3$ where we consistently draw high outliers when sampling $F_n$. This experiment brings the algorithm far from the theoretical requirements H1-H5 that assume $\nu > 12$.

\begin{figure}[]
\centering
\includegraphics[width=0.54\linewidth]{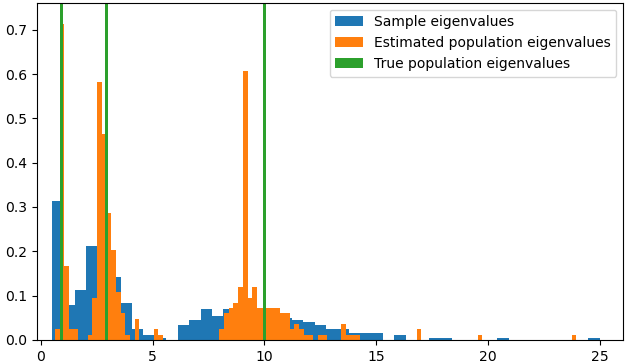}
\includegraphics[width=0.54\linewidth]{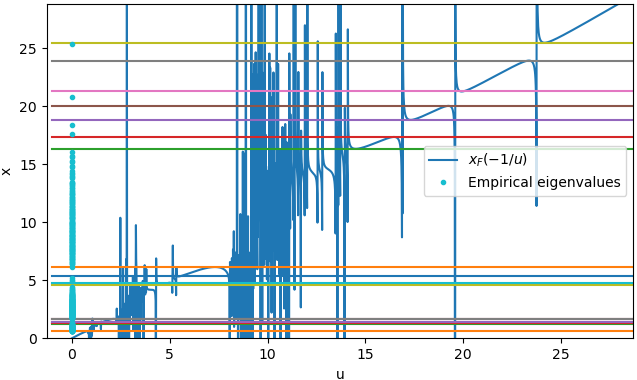}
\includegraphics[width=0.54\linewidth]{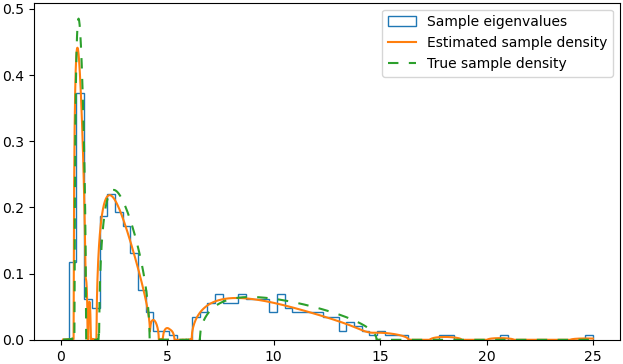}
\caption[Histograms, support identification and densities resulted from \textit{WeSpeR} (MD) procedure, Student noise.]{(Top) Histograms of sample eigenvalues, estimated population eigenvalues $\hat H$ with \textit{WeSpeR} (MD) procedure, and true population eigenvalues $H$. (Bottom) Estimated and true sample density computed on $S_F$ and sample eigenvalues' histogram. $H = \frac{1}{5} \mathbf{1}_{[1,\infty[} + \frac{2}{5} \mathbf{1}_{[3,\infty[} + \frac{2}{5} \mathbf{1}_{[10,\infty[}$, $D$ exponentially weighted with $\alpha = 1$, $c=0.1$, and heavy tails $Z_{ij} \sim t_3(0,1)$.}
\label{fig:wesper2}
\end{figure}

As shown in the Figure \ref{fig:wesper2}, these outliers in the sample spectrum affect badly the estimation of $H$, skewing it towards high values. Fortunately, as their frequency is quite low, only the higher values in $H$ are deteriorated and the estimated density $F'$ is still accurate in the core of the distribution. Some artifacts appear in the tail of the estimated $F'$ and $H$, with unwanted and isolated high values that fit the observed outliers. These outliers create small intervals in the estimated support $S_F$ around the observed outliers. 

In this extreme setting with $\nu=3$, far from the theoretical requirement, if one uses heavy tails, we recommend transforming the observed eigenvalues in order to reject the highest quantiles.

\section*{Acknowledgment}
We would like to thank Alexandre Miot and Gabriel Turinici for their insights and advices all along the work.

\bibliography{aistats_2025}

\newpage

\end{document}